\documentclass[a4paper, twoside]{amsart}

\usepackage{graphicx}
\usepackage{framed}

\usepackage[T1]{fontenc}							
\usepackage[latin1]{inputenc}

\usepackage[all]{xy}
\usepackage{pb-diagram, pb-xy}

\usepackage{amssymb}
\usepackage{amsthm}

\numberwithin{equation}{subsection}

\newtheoremstyle{myconj}
{0.8em}
{0.4em}
{\itshape}
{}
{\bfseries}
{.}
{ }
{\bfseries \thmname{#1} \thmnumber{#2}\thmnote{(#3)}}%

\newcounter{conjcount}

\theoremstyle{myconj}
\newtheorem{conj}[conjcount]{Conjecture}

\theoremstyle{plain}
\newtheorem*{thm*}{Theorem}
\newtheorem*{prop*}{Proposition}
\newtheorem*{lem*}{Lemma}

\newtheorem*{rem*}{Remark}
\newtheorem*{cor*}{Corollary}
\newtheorem*{ex*}{Example}
\newtheorem*{riglem*}{Rigidity lemma}
\newtheorem*{keylem*}{Key lemma}

\newcommand{\medrowheight}{\rule[-0.4em]{0em}{1.4em}}

\newcommand{\bbC}{{\mathbb C}}
\newcommand{\bbG}{{\mathbb G}}
\newcommand{\bbH}{{\mathbb H}}
\newcommand{\bbP}{{\mathbb P}}
\newcommand{\bbQ}{{\mathbb Q}}
\newcommand{\bbR}{{\mathbb R}}
\newcommand{\bbS}{{\mathbb S}}
\newcommand{\bbZ}{{\mathbb Z}}
\newcommand{\bbN}{{\mathbb N}}

\newcommand{\one}{{\mathbf{1}}}

\newcommand{\IC}{{\mathrm{IC}}}

\newcommand{\Qbar}{\overline{\mathbb Q}}

\newcommand{\calA}{{\mathcal A}}
\newcommand{\calB}{{\mathcal B}}
\newcommand{\calC}{{\mathcal C}}

\newcommand{\calG}{{\mathcal G}}
\newcommand{\calH}{{\mathcal H}}
\newcommand{\calJ}{{\mathcal J}}

\newcommand{\calM}{{\mathcal M}}
\newcommand{\calN}{{\mathcal N}}

\newcommand{\calT}{{\mathcal T}}
\newcommand{\calU}{{\mathcal U}}
\newcommand{\calV}{{\mathcal V}}
\newcommand{\calX}{{\mathcal X}}

\newcommand{\calZ}{{\mathcal Z}}
\newcommand{\calTheta}{{\Theta_{\calX}}}

\newcommand{\Pbar}{{\overline{P}}}

\newcommand{\spbar}{{\overline{sp}}}

\newcommand{\Ytilde}{{\tilde{Y}}}

\newcommand{\frakg}{{\mathfrak g}}

\newcommand{\Sp}{\mathrm{Sp}}
\newcommand{\Sl}{\mathrm{Sl}}

\newcommand{\SO}{\mathrm{SO}}
\newcommand{\MT}{\mathrm{MT}}
\newcommand{\MHM}{\mathrm{MHM}}

\newcommand{\Dbc}{{{D_{\hspace{-0.05em}c}}^{\hspace{-0.27em}b}}}
\newcommand{\Dbcbar}{\overline{D}{}^{\hspace{0.05em}b}_{\hspace{-0.05em}c}}
\newcommand{\Perv}{\mathrm{Perv}}
\newcommand{\pD}{{^p \! D}}
\newcommand{\stD}{{^{st} \! D}}
\newcommand{\ptau}{{^{p} \! \tau}}
\newcommand{\sttau}{{^{st} \! \tau}}
\newcommand{\pH}{{^p \! H}}

\newcommand{\pdelt}[1]{{^p\hspace{-0.05em} \delta_{#1}}}
\newcommand{\cdelt}[1]{{^c\hspace{-0.05em} \delta_{#1}}}

\newcommand{\End}{{\mathit{End}}}
\newcommand{\Hom}{{\mathit{Hom}}}
\newcommand{\Ext}{{\mathit{Ext}}}
\newcommand{\id}{{\mathit{id}\hspace{-0.1em}}}
\newcommand{\Rep}{{\mathrm{Rep}}}
\newcommand{\sRep}{{\mathrm{sRep}}}
\newcommand{\sgn}{{\mathit{sgn}}}

\begin{document}

\title[The symmetric square of the theta divisor]{The symmetric square of the theta divisor in genus 4}
\author{T. Kr\"amer and R. Weissauer}

\maketitle

\begin{abstract}
We study the variation of Hodge structures which arises from the intersections of two translates of theta divisors on a principally polarized abelian variety of dimension $g=4$, using a Tannakian description for the convolution product of perverse sheaves on abelian varieties.
\end{abstract}

\thispagestyle{empty}

\section*{Introduction}

Let $X$ be a complex principally polarized abelian variety (ppav for short) of dimension $g \geq 2$ with a symmetric theta divisor $\Theta \subset X$. It is well known that the geometry of the intersections 
\[
 Y_{x} \;=\; \Theta_{x} \cap \Theta_{-x},
\]
where $\Theta_x = \Theta + x$  denotes the translate of $\Theta \subset X$ by a point $x\in X(\bbC)$,
is closely related with Torelli's theorem \cite{DebTorelli} and with the Schottky problem~\cite{DeUpdate}. Here 
we are interested in the variation of \mbox{$\bbQ$-Hodge} structures given by $H^\bullet(Y_x, \bbQ)$ for varying~$x$. If the theta divisor is smooth, then for general $x$ the intersection $Y_x$ is a smooth variety of dimension $g-2$, and by the weak Lefschetz theorem the cohomology of $Y_x$ comes by restriction from the one of $X$ except in degree $g-2$.~So we mainly consider the quotient
\[
 H \;=\; H^{g-2}(Y_x, \bbQ) \, / \, H^{g-2}(X, \bbQ).
\]
The involution $\sigma=-id_X$ acts on~$Y_x$ and induces a decomposition $H = H_+ \oplus H_-$ into eigenspaces $H_\pm$ which are the fibres of two interesting variations $\calV_\pm$ of \mbox{$\bbQ$-Hodge} structures. These variations are studied in section~\ref{sec:vhs}, with a special emphasis on the case of low genus $g$. The low-dimensional cases in particular motivate the 

\begin{conj} \label{conj:vhs}
If $\Theta$ smooth, then $\calV_\pm$ are simple.
\end{conj}

Considered as Hodge modules in the sense of~\cite{SaIntro}, the variations $\calV_\pm$ of Hodge structures have two underlying perverse sheaves $\delta_\pm $ on $X$ which we introduce in section~\ref{sec:multiplier}. Thus conjecture~\ref{conj:vhs} would follow from

\begin{conj} \label{conj:thetasquare} 
If $\Theta$ is smooth, the perverse sheaves $\delta_\pm$ are simple.
\end{conj}

The main motivation for this conjecture comes from representation theory. Up to a skyscraper sheaf we construct $\delta_\pm$ in section~\ref{subsec:theta_square} as the symmetric respectively alternating square of the perverse intersection cohomology sheaf~$\delta_\Theta$ in the tensor category introduced in~\cite{WeT}. This category is equivalent to the category $\Rep(G)$ of algebraic representations of some (albeit in general unknown) complex algebraic group $G=G(X, \Theta)$, so conjecture~\ref{conj:thetasquare} would follow from

\begin{conj} \label{conj:thetagroup}
If $\Theta$ is smooth, then
\[
 G(X,\Theta) \;=\; 
 \begin{cases}
  \; \mathrm{SO}(g!, \bbC) & \textnormal{\em for $g$ odd}, \\
  \; \Sp(g!, \bbC) & \textnormal{\em for $g$ even},
 \end{cases}
\]
and $\delta_\Theta$ corresponds to the standard representation of this group.
\end{conj}

In fact, it is not hard to show that the groups on the right hand side give an upper bound for 
the group $G(X,\Theta)$; see section~\ref{subsec:thetagroupmotivation}.
In section~\ref{sec:super-mackey} we show that the two conjectures~\ref{conj:thetasquare} and~\ref{conj:thetagroup} are equivalent, using a classification of low-dimensional representations~\cite{KrW} and a Mackey argument. 

\medskip

Since every ppav of dimension $g\leq 3$ with a smooth theta divisor is a Jacobian variety, the above conjectures are true in these cases~\cite{WeBN}.  The first new case appears for $g=4$. If we drop the assumption that the theta divisor is smooth, this is also the first non-trivial case for the Schottky problem to characterize the Jacobian varieties among all ppav's.  The goal of the present paper\footnote{Recently we have been able by a different method to prove the main theorem for all $g$, see~\cite{KrW3} where we also discuss the relationship between conjecture~C and the Schottky problem.} is the

\medskip

{\bf Main theorem.} {\em For a general complex ppav $X$ of dimension $g\leq 4$ with a symmetric theta divisor $\Theta \subset X$ we have
\[
 G(X,\Theta) \;=\; 
 \begin{cases}
  \; \mathrm{SO}(g!, \bbC) & \textnormal{\em for $g$ odd}, \\
  \; \Sp(g!, \bbC) & \textnormal{\em for $g$ even},
 \end{cases}
\]
and $\delta_\Theta$ corresponds to the standard representation of this group.}

\medskip

Here by the term {\em general} we mean that the claim holds for every ppav in a suitable Zariski-open dense subset of the moduli space $\calA_g$ of complex ppav's of dimension~$g$. The proof of the main theorem will be given in sections~\ref{sec:outline} and \ref{sec:decomposition} and uses a degeneration of $X$ into the Jacobian of a general curve. Such a degeneration gives a restriction functor $\rho: \Rep(G(X, \Theta)) \to \Rep(G_\Psi)$, where $G_\Psi\subset G(X, \Theta)$ is an algebraic subgroup defined via the formalism of the nearby cycles $\Psi$ as outlined in section~\ref{sec:cycles}. Together with the information obtained from this restriction functor we consider a second tensor functor $\MT$ from $\Rep(G(X, \Theta))$ to the category of representations of the Mumford-Tate group of the theta divisor. To make effective use of this second functor, we exploit the fact that the primitive cohomology of the theta divisor is related to the intermediate Jacobian of a cubic threefold ~\cite{Do}, \cite{Iz} and use results of~\cite{Co}, \cite{CM} and \cite{IvS} on cubic threefolds, see section~\ref{sec:motivic}.

\section{Variations of Hodge structures (conjecture~\ref{conj:vhs})} \label{sec:vhs}

In this section we assume that $\Theta$ is smooth. We study the variations~$\calV_\pm$ of Hodge structures and their Hodge decomposition for $g\leq 4$.

\subsection{Lemma} \label{subsec:bertini}
{\em For generic $x\in X(\bbC)$
the translates $\Theta_{x}$ and $\Theta_{-x}$ intersect each other transversely, hence $Y_x$ is smooth.}

\medskip

{\em Proof.} $\Theta$ is defined by the zero locus of the Riemann theta function $\theta(z)=\theta(\tau, z)$ on the 
universal covering $p: \bbC^g\rightarrow \bbC^g/(\bbZ^g+\tau \bbZ^g) \cong X(\bbC)$. For smooth $\Theta$ the gradient $\theta'(z)$ is non-zero for all $z\in \bbC^g$ with $\theta(z)=0$. If our claim were false, then we could find a non-empty analytic open subset $V \subset \bbC^g$ and a complex analytic map $s: V \rightarrow \bbC^g$ such that $\Theta_{p(z)}$ and $\Theta_{-p(z)}$ intersect non-transversely in~$s(z)$ for all $z\in V$, i.e.~such that
\begin{gather} \label{eq:onZ}
\qquad \theta(s(z) + z)  \;=\; \theta(s(z) - z) \;=\; 0, \\
\label{eq:tangential}
\theta'(s(z)+ z) \;=\; \lambda(z)\cdot \theta'(s(z)-z)
\quad
\textnormal{for some $\lambda(z)\in \bbC^*$},
\end{gather}
for all $z\in V$. With the $g\times g$ unit matrix $E$ and  the Jacobian matrix $(Ds)(z)$ of~$s$ at $z$, taking gradients of~(\ref{eq:onZ}) implies
\begin{gather} \label{eq:plus}
 (E+(Ds)(z))\cdot \theta'(s(z)+z) \;=\; 0, \\
\label{eq:minus}
 (E-(Ds)(z))\cdot \theta'(s(z)-z) \;=\; 0.
\end{gather}
If we multiply equation (\ref{eq:minus}) by $\lambda(z)$, plug in (\ref{eq:tangential}) and add (\ref{eq:plus}), we obtain that $\theta'(s(z)+z)=0$ which  contradicts the smoothness of $\Theta$. \qed

\subsection{Variations of Hodge structures} To construct the variations of Hodge structures from the introduction, we realize the intersections $Y_x=\Theta_x \cap \Theta_{-x}$ as the fibres of a smooth proper family.

\begin{lem*} \label{subsec:family}
Over some Zariski-open dense subset $U\subset X$ there exists a smooth proper family $\pi: Y_U \to U$ with fibres $\pi^{-1}(2x) \cong Y_x$ and an \'etale involution $\sigma: Y_U\to Y_U$ with $\sigma|_{\pi^{-1}(2x)} \cong (-\id_X)|_{Y_x}$ for $2x$ in $U(\bbC)$.
\end{lem*}

\medskip

{\em Proof.} Let $\pi: \Theta\times \Theta \rightarrow X$ be the addition. Projecting from $\pi^{-1}(2x)\subset \Theta \times \Theta$ onto the first factor and translating by $x$, we get an isomorphism $\varphi: \pi^{-1}(2x) \cong  \Theta_x\cap \Theta_{-x}$ under which the morphism $\sigma: \Theta\times \Theta \rightarrow \Theta\times \Theta, \, (t_1, t_2) \mapsto (t_2,t_1)$ becomes identified with the involution $(-\id_X)|_{\Theta_x\cap \Theta_{-x}}$. By section~\ref{subsec:bertini} the fibres of $\pi$ are generically smooth. So there exists an open dense subset $U\subset X$ such that for $Y_U = \pi^{-1}(U)$ the restriction $\pi = \pi|_{Y_U}$ is smooth, and for sufficiently small $U$ the involution $\sigma$ will be \'etale on~$Y_U$. \qed

\medskip

\label{subsec:vhs}
For all $\nu \in \bbZ$, the higher direct images
$R^\nu\pi_*(\bbQ_{Y_U})$ are variations of $\bbQ$-Hodge structures \cite[cor.~10.32]{PS}.
There exists a constant subvariation 
$
 H^\nu \hookrightarrow R^\nu \pi_*(\bbQ_{Y_U}) 
$
such that for $2x\in U(\bbC)$, the fibre $H^\nu \subset R^\nu \pi_*(\bbQ_{Y_U})_{2x} = H^\nu(Y_x, \bbQ)$ is the pull-back of the cohomology of $X$ to~$Y_x$. In particular $H^\nu \cong H^\nu(Y_x,\bbQ)$ for all $\nu\neq g-2$. In what follows we define
\[ \calV = R^{g-2}\pi_*(\bbQ_{Y_U}) / H^{g-2}
\]
and consider the eigenspace of $\sigma^*$
\[
 \calV_\pm = \ker(\sigma^* \mp (-1)^{g-1} \id_\calV)
\]
with respect to the eigenvalues~$\pm (-1)^{g-1}$. The reason for this choice of signs will become clear in section~\ref{subsec:theta_square}.

\subsection{Small-dimensional cases} \label{subsec:g3}
In this section we fix $x\in X(\bbC)$ with $2x\in U(\bbC)$. For
$Y = Y_x$, we define $Y^+$ as the quotient $Y/\langle \sigma \rangle$ of $Y$ by the involution $\sigma=-\id_X|_Y$.

\medskip

{\em The case $g=2$}.  Here $Y$ consists of $2$ points, $Y^+$ is a single point. Hence ${\mathcal V}_-=0$ and~${\mathcal V}_+=\bbQ_U$ is the constant variation with Hodge degree $(0,0)$. 

\medskip

{\em The case $g=3$}. Here~$(X, \Theta)$ is the Jacobian of a smooth curve $C$ and~$Y^+\to Y$ is an \'etale double cover of smooth curves with $Y$ of genus~$7$. To any \'etale double cover of curves one may associate a ppav, its Prym variety \cite{MuPrym} \cite[ch.~12]{BL}.
It turns out that for generic $(X, \Theta)$ the Prym variety~$P$ of the cover $Y\rightarrow Y^+$ is isomorphic to~$(X,\Theta)$. Every \'etale double cover with this Prym variety arises like this for some~$x$. Furthermore, the coverings for two points $x_1, x_2$ are isomorphic iff $x_1 = \pm x_2$. Hence the \'etale double covers with given Prym variety $(X,\Theta)$ are parametrized by an open dense subset $W$ of the Kummer variety $X/\langle \pm 1\rangle$. Points outside $W$ parametrize degenerate double covers.
By construction of Prym varieties, the Jacobian $JY$ is isogenous to the product $P\times JY^+$, hence
$
 H^1(Y, \bbC) \cong H^1(X,\bbC)\oplus H^1(Y^+,\bbC).
$
So it follows that~${\mathcal V}_-=0$, and ${\mathcal V}_+$ is a variation of Hodge structures of abelian type whose fibres $H^1(Y^+, \bbQ)=\bbQ^8$ have Hodge numbers $h^{1,0}=h^{0,1}=4$. 
This construction exhibits the Prym locus in $\calA_4$ birationally as the universal Kummer variety fibered over the 
moduli space $\calA_3$, studied in~\cite{Re}, \cite{Kr1}.
\medskip

\label{subsec:g4}
{\em The case $g=4$}. Here $Y\rightarrow Y^+$ is an \'etale covering of smooth surfaces. 
Via the Gauss-Bonnet and Grothendieck-Riemann-Roch formulae one checks that $\calV_+$ has rank $52$, with Hodge numbers $h^{2,0}=h^{0,2}=11$ and $h^{1,1}=30$, whereas~$\calV_-$ has rank $6$ and pure Hodge type~$(1,1)$. Again~$\calV_-$ seems to be the more accessible one of the two variations of Hodge structures. 
Indeed, by the Lefschetz $(1,1)$-theorem the fibers of~$\calV_-$ are spanned by six cycles on $Y^+$ generating $H^{1,1}(Y^+, \bbC)/H^{1,1}(X, \bbC)$.
So~$\calV_-$ has an underlying finite monodromy group $\Gamma=\Gamma(X, \Theta)$ with a six-dimensional faithful representation. It turns out that $\Gamma$ is the Weyl group $W(E_6)$ or its finite simple subgroup of index~$2$. 
To show this one uses the classification of lattices with small discriminant to see that the fibres of~$\calV_-$ have the underlying N\'eron-Severi lattice $E_6(-1)$, so $\Gamma$ is a subgroup of $Aut(E_6)=W(E_6)\times \{ \pm 1\}$. From the intersection configuration of the $27$ Prym-embedded curves in $Y$ of~\cite[sect.~4.3]{Iz}, one furthermore deduces that the projection $W(E_6)\times \{\pm 1\} \to W(E_6)$ maps $\Gamma$ isomorphically onto a subgroup of $W(E_6)$. On the other hand, one obtains a lower bound on $\Gamma$ by a degeneration of $(X, \Theta)$ into a Jacobian variety $(X_0, \Theta_0)$. The monodromy group $\Gamma_0$, underlying the analog of $\calV_-$ on $(X_0, \Theta_0)$, is a subquotient of~$\Gamma$. In the Jacobian case the configuration of the $27$ Prym curves of loc.~cit.~is no longer symmetric; precisely $12$ of them are smooth. They come in $6$ pairs of curves which are interchanged by the involution $\pm \id_{X_0}$, and the associated $6$ cycles are permuted by the monodromy operation of $\Gamma_0$. From this one deduces that $\Gamma_0$ contains the alternating group $A_5$. Altogether this already forces $\Gamma$ to be the Weyl group $W(E_6)$ or its finite simple subgroup of index~$2$, see~\cite{Kr2}.

\section{The Mumford-Tate group $\MT(\Theta)$ in Genus 4} \label{sec:motivic}

One of the ingredients to the proof of the main theorem will be the analysis of the Hodge structure on tensor products of $\bbH^\bullet(X, \delta_\Theta)$. As we recall below, this Hodge structure is controlled by the corresponding Mumford-Tate group $\MT(\Theta)$. 
For reductive groups $G$ let $G_{sc}$ be the simply connected covering of the derived group $[G^0, G^0]$ of its Zariski connected component~$G^0$. If a subset of an algebraic variety is contained in a countable union of proper closed subvarieties, it is said to be a {\em meager} subset. The goal of this section is to prove the following result.

\subsection{Theorem.}  \label{thm:MTB}
{\em For every ppav $(X, \Theta)$ outside some meager subset of $\calA_4$ we have the Mumford-Tate groups
\[ \MT(\Theta)_{sc} = \MT(X) \times \Sp(10, \bbC)
\quad
\textnormal{\em and}
\quad
\MT(X)=\Sp(8, \bbC), \]
so that we have canonical homomorphisms
\[
 \xymatrix@=2em{
 \MT(X) \ar[r] \ar[dr] & \MT(\Theta)_{sc} \ar[d] & \Sp(10, \bbC) \ar[l] \ar[dl] \\
 & \MT(\Theta) &
}
\]
If we view $V=\bbH^\bullet(X, \delta_\Theta)$ as a representation of $\MT(\Theta)$, the action of $\MT(X)$ on~$V$ describes the part of the cohomology which comes by restriction from $X$. The action of $\Sp(10, \bbC)$ describes the remaining part.
}

\medskip

Recall from~\cite[sect.~2.1]{DelHII} that giving a \mbox{$\bbQ$-Hodge} structure is tantamount to giving a finite-dimensional vector space $V$ over $\bbQ$ together with a homomorphism $h: \bbS = \mathrm{Res}_{\bbC/\bbR}(\bbG_{m, \bbC}) \rightarrow \mathrm{Gl}(V)_\bbR$ of real algebraic groups whose composite with the weight cocharacter of $\bbS$ is defined over $\bbQ$.
The Mumford-Tate group $\MT(V)$ is the smallest algebraic subgroup of $\mathrm{Gl}(V)$ over~$\bbQ$ over which $h$ factors. The $\bbQ$-Hodge substructures of any tensor power of $V$ are precisely its $\MT(V)$-stable subspaces. Note that the Mumford-Tate group of polarized $\bbQ$-Hodge structures is reductive~\cite[prop.~3.6]{DelH}; this in particular applies to $\MT(X):=\MT(H^\bullet(X, \bbQ))$ for smooth projective varieties $X$ over~$\bbC$.

\subsection{The primitive part of the cohomology} \label{subsec:theta-cohomology}
Let $\Lambda^i$ denote the $i$-th exterior power of the standard representation of $\Sp(2g, \bbQ)$. Then, by~\cite[prop.~17.3.2]{BL}
for every ppav $(X, \Theta)$ outside a meager subset of~$\calA_g$ we know  that
\[ \MT(X) \;=\; \MT(X)_{sc} \;=\; \Sp(2g, \bbQ)
 \quad
 \textnormal{and}
  \quad
  H^i(X, \bbQ)=\Lambda^i \ . \]
If $\Theta$ is smooth, the weak Lefschetz theorem shows that for $0\leq |i| \leq g-1$ 
and some $\bbQ$-Hodge structure~$B$
\[
 H^{g-1+i}(\Theta, \bbQ) \;=\;
 \begin{cases}
 \; \Lambda^{g-1-|i|} & \textnormal{for $i\neq 0$}, \\
 \; \Lambda^{g-1} \, \oplus \, B & \textnormal{for $i=0$}\ .
 \end{cases}
\]
Hence
$\MT(\Theta) = \MT\bigl(H^\bullet(X, \bbQ) \oplus B \bigr)$.
For any $\bbQ$-Hodge structures~$V_1$, $V_2$ and $V=V_1\oplus V_2$ we have a closed embedding
$\iota: \MT(V) \hookrightarrow \MT(V_1) \times \MT(V_2)$, and
the image of $\iota$ surjects onto each of the two factors. If the Mumford-Tate groups are reductive, this surjectivity carries over to $\MT(-)_{sc}$, and the kernel of the induced map $\iota_{sc}: \MT(V)_{sc} \to \MT(V_1)_{sc} \times \MT(V_2)_{sc}$ is contained in the kernel of the natural map $\MT(V)_{sc} \to \MT(V)$.
Hence we have a commutative diagram
\[
\xymatrix@M=0.5em{
 & \MT(\Theta)_{sc} \ar[d]^-{\iota_{sc}} \ar@{->>}[dl] \ar@{->>}[dr] & \\
 \Sp(2g, \bbQ) & \Sp(2g, \bbQ) \times \MT(B)_{sc} \ar@{->>}[l]_-{p_1} \ar@{->>}[r]^-{p_2} & \MT(B)_{sc} 
}
\]
with $\ker(\iota_{sc})$ contained in $\ker(\MT(\Theta)_{sc} \to \MT(\Theta))$.

\begin{lem*}
If for any ppav $X$ outside a meager subset of $\calA_4$ the Hodge structure~$B$ satisfies $\End_{\MT(B)}=\bbQ$, then theorem~\ref{thm:MTB} holds.
\end{lem*}

{\em Proof.}  Since $B$ is a Hodge structure of abelian type, by~\cite[th.~1]{Ri} the condition $\End_{\MT(B)}(B)=\bbQ$ implies $\MT(B)=\MT(B)_{sc}=\Sp(10, \bbC)$. Since the projection $\MT(\Theta)_{sc} \to \MT(X)=\Sp(8, \bbC)$ is surjective, there exists a reductive group~$G$ with
$\MT(\Theta)_{sc} = \Sp(8, \bbC) \times G$. Hence the simplicity of the Hodge structure $B$ forces  
$B\cong B_1 \boxtimes B_2$
with irreducible representations $B_1$ of $\Sp(8, \bbC)$ and $B_2$ of~$G$. If~$B_1$ were non-trivial, then $\dim(B)=10$ would imply that $\dim(B_1) \in \{ 2, 5, 10\}$ which is impossible~\cite{AEV}. Hence $B_1$ is the trivial representation, $B_2$ is the standard representation of $\MT(B)=\Sp(10, \bbC)$, and $\iota_{sc}$ is an isomorphism. \qed \\

For every ppav $(X, \Theta)$ outside a meager subset of $\calA_4$, the Hodge structure $B$ indeed satisfies $\End_{\MT(B)}(B)=\bbQ$. This is shown in the remaining part of this section, and by the lemma above completes the computation the Mumford-Tate group. Notice that $\MT(B)\subseteq \Sp(10, \bbC)$ by ~\cite[B.62]{Go}. Since the Mumford-Tate group of a variation of Hodge structures is constant outside the complement of some meager subset and only becomes smaller on this meager subset, it suffices to prove $\End_{\MT(B)}(B) = \bbQ$ for a single ppav $B$. For this we consider intermediate Jacobians of cubic threefolds.

\subsection{Intermediate Jacobians of cubic threefolds} \label{subsec:intermediate-jacobians}
Recall from~\cite{CG} that for any smooth cubic threefold $T\subset \bbP^4_\bbC$ the intermediate Jacobian
\[
 JT \;=\; H^{2,1}(T)^* / H_3(T, \bbZ)
\]
is a simple ppav of dimension $5$ and determines $T$ up to isomorphism. 
Let us denote by $\calT$ the closure  in the moduli space $\calA_5$ of the locus of all these intermediate Jacobians.
In~\cite{Do},
to each ppav $(X, \Theta)$ in some Zariski-open dense subset~$\calA_4^\circ$ of~$\calA_4$ a smooth cubic threefold $T=T_{(X, \Theta)}$ has been associated together with a $2$-torsion point $\mu$ of its intermediate Jacobian $JT$. The associated map
$\varphi:  \calA_4^\circ  \rightarrow  \calT$
is generically finite and dominant by \cite[birationality of $\chi$ in thm.~5.2 (2)]{Do}. Hence the image of $\varphi$ contains a Zariski-open dense subset of $\calT$. 
For this construction Donagi analyzes the Prym map,
associating to curves of genus 5 together with a two-torsion point in their Jacobian
a Prym variety in $\calA_4$, and the fibers $\tilde F$ of this Prym map over given points $(X,\Theta)$ in $\calA_4^0$. 
By studying Prym varieties for plane quintic curves of genus six, Donagi defines an involution on the Prym fibers $\tilde F$. He shows that the associated quotient
map $\tilde F \to F$ is an etale double covering. The quotient $F$ turns out to be the Fano surface embedded into $JT$, parametrizing the projective lines in the cubic threefold $T= \varphi(X,\Theta)$. In fact $Alb(F) \cong JT$, whereas the isogenous abelian variety $Alb(\tilde F)$ is isomorphic to $B$ via a higher
Abel-Jacobi map defined by a family $\calC(X,\Theta) \to \tilde F$ of stable curves in $\Theta$ as defined in \cite{Iz}, p.133, parametrized by the base $\tilde F$. 
Using a degeneration argument one can show that this Abel-Jacobi map is nontrivial.
If on the other hand $JT$ is a simple abelian variety, this map is an isogeny
(and then even an isomorphism). Replacing the analytic definition of~$JT$ from~\cite{CG} by an algebraic definition, viewing $JT$ as quotient of the Chow group $A(T)$, the above constructions can be done over any algebraically closed field of characteristic zero. Hence if $(X,\Theta)$ is defined over $\overline\bbQ$, then so is $T$ and the isogeny between $B$ and $JT$.

\label{subsec:End0_and_simplicity}

\medskip
For the moment we only use that for $(X, \Theta)\in \calA_4^\circ(\bbC)$ the $\bbQ$-Hodge structures on $H^1(B,\bbQ)$
and $H^1(JT_{(X, \Theta)}, \bbQ)(-1)$ are isomorphic, as shown in~\cite{IvS} and~\cite{Iz}. 

\medskip
For abelian varieties $A$ put $\End^0(A)=\End(A)\otimes_\bbZ \bbQ$.  Then a basic property of Mumford-Tate groups~\cite[B.60]{Go} shows
$
 \End_{\MT(B)}(B) \;=\; \End^0(JT_{(X, \Theta)}).
$
Thus to verify the condition of the last lemma, it remains to show $$\End^0(JT_{(X, \Theta)}) = \bbQ\ $$ for a suitable $(X, \Theta)$ in $\calA_4^\circ$. By the above
discussion this becomes a question on the five-dimensional intermediate Jacobians of generic cubic threefolds. So we can now use a degeneration argument of Collino for cubic threefolds, which in turn reduces our task to the study of a certain extension class of a generic Jacobian variety $JC$ of genus 4 by a torus.

\subsection{Collino's family}
By~\cite[part II]{Co} and~\cite[p.~44-45]{CM} there exists a group scheme 
$ \calB \rightarrow S $
over a non-empty Zariski-open subset $S\subset \bbP^1_\bbC$ and a point $s\in S(\bbC)$ such that
\begin{itemize}
 \item for all $t \neq s$ one has $\calB_t \cong JT_{(X_t, \Theta_t)}$ for some $(X_t, \Theta_t) \in \calA_4^\circ$, 
 \item the special fibre $\calB_s$ is an extension
\[
 E: \; \; 0 \longrightarrow \bbG_m \longrightarrow \calB_s \longrightarrow JC \longrightarrow 0
\]
for some general curve $C$ of genus $4$. The theta divisor $JC$ has precisely two singular points~$\pm e$. The class of the extension $E$ in $\Ext(JC, \bbG_m)$ is mapped under the isomorphism $\Ext(JC, \bbG_m) \stackrel{\sim}{\longrightarrow} \mathrm{Pic}^0(JC) \stackrel{\sim}{\longrightarrow} JC$ from \cite[16.VII]{Se} to the point $2e$ or $-2e$. Since $C$ is general, we know that $e\neq -e$ and the extension $E$ is non-trivial.
\end{itemize}
We will show $\End^0(\calB_t) = \bbQ$ for all $t$ outside a meager subset of $S(\bbC)$. Note that over $S^*=S\setminus \{s\}$ the family $\calB\to S$ restricts to an abelian scheme $\calB^*\to S^*$. So by \cite[prop.~7.5]{DelK3} and \cite[4.1.3.2]{DelHII} the restriction map $\End^0(\calB^*/S^*)\to \End^0(\calB_t)$ is an isomorphism for all but {\em countably} many $t\in S^*(\bbC)$ (this is why we had to exclude a meager subset of ppav's in the theorem). So it suffices to show $\End(\calB^*/S^*)=\bbZ$, and this is equivalent to the claim in~\ref{subsec:endo-generic} below.

\subsection{Endomorphisms of the special fibre} \label{subsec:End_at_0}
We first show that for general choice of $C$ one has
\[ \End(\calB_s) = \bbZ. \]

{\em Proof.} Every $\psi\in \End(\calB_s)$ preserves the toric part $\bbG_m\subset \calB_s$ and induces an endomorphism $\psi_{JC}$ of $JC = \calB_s/\bbG_m$. So we have a ring homomorphism
\[
 (-)_{JC}: \; \End(\calB_s) \; \longrightarrow \; \End(JC), \quad \psi \mapsto \psi_{JC}.
\]
For general $C$ we know $\End(JC)=\bbZ$, and then $(-)_{JC}$ is surjective because its image contains $1=(\id_{\calB_s})_{JC}$. Now suppose $\psi \in \End(\calB_s)$ and $\psi_{JC}=0$, i.e.~$\psi$ factors over $\bbG_m \subset \calB_s$. Then $\psi|_{\bbG_m}$ is a character $z \mapsto z^n$ of $\bbG_m$ for some $n\in \bbZ$. If $\psi \neq 0$, we must have $n\neq 0$. However, the image of the restriction map
\[
 \mathit{res}: \; \Hom(\calB_s, \bbG_m) \; \longrightarrow \; \Hom(\bbG_m, \bbG_m) \; = \; \bbZ
\]
contains $n = \mathrm{res}(\psi)$, so the image of the differential $d$ of the $\Ext$-sequence
\[
 \cdots \stackrel{\mathit{res}}{\longrightarrow} \; \bbZ \;=\; \Hom(\bbG_m, \bbG_m) \; \stackrel{d}{\longrightarrow} \; \Ext(JC, \bbG_m) \; \longrightarrow \; \cdots
\]
is a quotient of $\bbZ/n\bbZ$. In particular, the image of $d$ is an $n$-torsion group. But by construction of the $\Ext$-sequence, $d(\id_{\bbG_m})$ is the class of the extension defining the semiabelian variety $\calB_s$. Hence $\pm 2 e \neq 0$ is an $n$-torsion point of $JC$. For $C$ varying this contradicts lemma~\ref{lem:torsion} below. Hence $(-)_{JC}$ is an isomorphism.
\qed

\subsection{Torsion points} \label{lem:torsion}
Let $\calM_4$ be the moduli space of smooth projective curves of genus $4$. Over some Zariski-open dense subset $\calU$ of $\calM_4$, the map $C\mapsto \pm 2e$ defines a section of the universal Jacobian variety, and this is not the zero section $0_{\,\calU}$.

\begin{lem*}
If $\calU \subset \calM_4$ is a Zariski-open dense subset and $\Sigma\neq 0_{\,\calU}$ is a section of the universal Jacobian variety $\pi: \calX \rightarrow \calU$, then~$\Sigma$ defines a non-torsion point in all fibres of $\pi$ over the complement of a meager subset of~$\calU$.
\end{lem*}

{\em Proof.} For $m\in \bbZ$, let $\calZ_m$ denote the zero locus of $m\cdot \Sigma$. Then $\calZ_m$ is a Zariski closed subset of~$\calU$, and we must show it is not all of $\calU$. If it were for some $m\geq 0$, then for minimal such $m$ the section $\Sigma$ would define a section over some Zariski open dense subset of $\calM_4$ to the map
$_m \calM_4 \longrightarrow \calM_4$,
where $_m\calM_4$ is the moduli space of cyclic \'etale covers of precise order $m$ of curves of genus $4$. Then $_m\calM_4$ would be reducible, contradicting~\cite{BF}.
\qed

\subsection{Endomorphisms of the generic fibre} \label{subsec:endo-generic}
For the generic point $\eta$ of $S$ we now claim that
\[ \End(\calB_\eta) = \bbZ. \]

{\em Proof.} Let $\calN$ be the N\'eron model $\calN$ of $\calB_\eta$ over $S$. Its universal property gives an $S$-morphism $\calB \rightarrow \calN$. Since this morphism induces an isomorphism of the generic fibre and since $\calB_{s_0}$ is semi-abelian, by~\cite[prop. 7.4.3]{BLR} it induces an isomorphism
of~$\calB$
onto the connected component $\calN^0 \subset \calN$. In section~\ref{subsec:End_at_0} we have shown $\End(\calB_{s_0})=\bbZ$. Now consider the composite ring homomorphism
\[
 \varphi: \;
 \End^0(\calB_\eta) \; \stackrel{\sim}{\longrightarrow} \; \End^0(\calN_\eta) 
 \; \stackrel{\sim}{\longrightarrow} \; \End^0(\calN/S) 
 \; \stackrel{\mathit{res}}{\longrightarrow} \; \End^0(\calB_s) \; = \; \bbQ\ ,
\]
where the first isomorphism comes from the identification $\calB_\eta \stackrel{\sim}{\longrightarrow} \calN_\eta$, the second isomorphism is due to the universal property of $\calN$, and the third map $\mathit{res}$ denotes restriction to the fibre $\calB_s = (\calN^0)_s$.
The image of $\varphi$ contains $1=\varphi(\id_{\calB_\eta})$, so $\varphi$ is surjective. 
On the other hand, $\calB_\eta$ is a simple abelian variety; indeed, all $\calB_t$ with $t\neq s$ are intermediate Jacobians of smooth cubic threefolds, hence simple by~\cite{CG}. Therefore $\End^0(\calB_\eta)$ is a skew field, and since $\varphi$ is a surjective ring homomorphism, it follows that $\ker(\varphi)=0$. Thus $\End^0(\calB_\eta)=\bbQ$, and our claim follows.
\qed \\

If combined, sections 2.4-2.7 show that $B$ is simple for generic $(X,\Theta)$
respectively that $JT$ is simple for a generic cubic threefold $T$. 
Hence the conditions of the last lemma have now been verified.

\subsection{$\ell$-adic counterpart} \label{motivic} Instead of the Mumford-Tate groups $MT(\Theta)$, $MT(X)$ one can consider the Zariski closure of the geometric \'etale fundamental group acting
on the $\overline \bbQ_{\ell}$-\'etale cohomology to define the motivic Galois groups 
$M(\Theta)$, $M(X)$. Then, one can similarly prove
$$ M(\Theta)_{sc} = M(X) \times Sp(10,\overline \bbQ_{\ell})\ \mbox{ and }
M(X)= Sp(8,\overline \bbQ_{\ell})\ $$
for a generic ppav $X$ of genus 4.  
To see this, notice that the image $G$
of $M(\Theta)_{sc}$ under the projection $M(X) \times Sp(10,\overline \bbQ_{\ell}) \to Sp(10,\overline \bbQ_{\ell})$ 
contains a subgroup which is isomorphic to~$Sp(8,\overline \bbQ_{\ell})$ (by a degeneration to the curve case). If $G \neq Sp(10,\overline \bbQ_{\ell})$, group-theoretic arguments would imply that the normalizer $N$ of the subgroup~$G$ in $Sp(10,\overline \bbQ_{\ell})$ is contained in a proper Levi subgroup $L$ of $Sp(10,\overline \bbQ_{\ell})$ which acts reducibly on the 10-dimensional fundamental representation of $Sp(10,\overline \bbQ_{\ell})$. But this is impossible, since by a specialization argument this would also hold for all points in a Zariski-dense open subset.
Consider such a point of $\calT$ defined over $\overline \bbQ$. It is already defined
over some number field, say $E$. Since the arithmetic $\ell$-adic monodromy group of $JT$ normalizes the geometric $\ell$-adic monodromy group of
$JT$, the arithmetic $\ell$-adic monodromy group of
$JT$ over $E$ must therefore be contained in $L(\overline \bbQ_{\ell})$. Hence by Falting's theorem $JT$ is not simple.  

\medskip
Since the points of $\calT$ defined over $\overline \bbQ$ are Zariski dense in $\calT$, 
this provides us with a Zariski-dense subset of points of $\calT$ for which $JT$ is not simple. 
For a generic threefold $T$ in $\calT$ the abelian variety $JT$ is simple by the conclusions of sections~\ref{subsec:intermediate-jacobians} through~\ref{subsec:endo-generic} above, and therefore the locus of decomposable $JT$ in $\calT$ can not be Zariski dense in $\calT$. This contradiction shows that $G=Sp(10,\overline \bbQ_{\ell})$, and the rest of the argument is analogous to the proof of the lemma in section 2.2.

\subsection{Higher genus} \label{subsec:motivic-curve-case} 
For the rest of this section we drop the assumption $g=4$. It seems likely that
for general $(X,\Theta)$ the cohomology 
$$H^{g-1}(\Theta,\bbQ)/H^{g-1}(X,\bbQ)$$  
decomposes into at least $\lfloor \frac{g-1}{2} \rfloor$ simple $\bbQ$-Hodge substructures. To illustrate this let us consider the case of Jacobians. Let $C$ be a general curve of genus $g\geq 3$. For $\nu \in \bbZ$ let $\Lambda^\nu$ be the $\nu$-th exterior power of the standard representation of $\MT(\Theta_{JC})_{sc} =\Sp(2g, \bbQ)$, and consider $\bbH^\nu(JC, \delta_{\Theta_{JC}})$ as a representation of the group $\langle \sigma^* \rangle \times \MT(\Theta_{JC})_{sc}$ for the involution $\sigma = -\id_{JC}$. Let $\sgn$ be the nontrivial character of $\langle \sigma^*\rangle$.

\begin{lem*}
For $|\nu| \leq g-1$ define $n(\nu) = g-1-|\nu|$. Then 
\[ 
 \bbH^\nu (JC, \delta_{\Theta_{JC}}) \ =\  \bigoplus_{\mu = 0}^{\lfloor \frac{n(\nu)}{2} \rfloor} \\sgn^{n(\nu)+\mu} \boxtimes \Lambda^{n(\nu)-2\mu}
  \ .
\]
\end{lem*}

{\em Proof.}
Up to the $\sigma$-action this is clear as $\bbH^\bullet(JC, \delta_{\Theta_{JC}}) = \Lambda^{g-1} (\bbH^\bullet(JC, \delta_C))$ by the non-hyperelliptic case of \cite[cor.~13(iii) on p.~64 and p.~124]{WeBN}, but note that the perverse sheaf $\delta_C$ is not $\sigma$-equivariant!
To find the $\sigma$-action on $\bbH^\bullet(JC, \delta_{\Theta_{JC}})$, we 
use the evaluation map \cite[sect.~4.2]{BrB}
\[ 
 H^\bullet (JC, \bbQ)\otimes_\bbQ \bbQ[x] 
 \;\longrightarrow\;
 H^\bullet(C^{(g-1)}, \bbQ) 
 \;\cong\;
 \bbH^{\bullet}(JC, \delta_{\Theta_{JC}})[1-g]
\]
which is an isomorphism in degrees $\leq g-1$ and $\sigma$-equivariant, where on the left hand side $\sigma$ acts in the usual way on $H^\bullet(JC, \bbQ)$ and on the powers of the variable~$x$ it acts by
\[
 \sigma^*(x^\mu) \;=\; \sum_{i=0}^\mu (-1)^i \frac{[\Theta_{JC}]^{\mu - i} \otimes x^i}{(\mu-i)!} \, 
\]
as observed in prop.~4.3.1 of loc.~cit. From this our claim easily follows.
\qed

\section{Convolutions of perverse sheaves (conjecture~\ref{conj:thetasquare})} \label{sec:multiplier}

In this section we introduce the perverse sheaves $\delta_\pm$ that occur in the formulation of conjecture~\ref{conj:thetasquare}.

\subsection{The convolution product} \label{sec:convolution}
Put $k=\bbC$ or $k=\Qbar_l$ for a fixed prime $l$, and denote by $\Dbc(X,k)$ the triangulated category of bounded constructible complexes of sheaves with coefficients in $k$ as in \cite{KW}. The group law $a:X\times X\rightarrow X$ induces a convolution product on this triangulated category via the formula
\[
 \gamma_1 * \gamma_2  =  Ra_*(\gamma_1\boxtimes \gamma_2)
 \quad \textnormal{for} \quad
 \gamma_1, \gamma_2 \in \Dbc(X, k).
\]
With respect to this convolution product, the category $\Dbc(X, k)$ becomes a \mbox{$k$-linear} tensor category~\cite[sect.~2.1]{WeBN} 
whose unit object is the skyscraper sheaf~$\delta_{\{0\}} = \one$ supported in the origin. For hypercohomology the relative K\"unneth formula implies $$\bbH^\bullet(X, \gamma_1 * \gamma_2) = \bbH^\bullet(X, \gamma_1) \otimes \bbH^\bullet(X, \gamma_2)\ .$$
Let $\Perv(X)\subset \Dbc(X, k)$ be the abelian category of perverse sheaves on~$X$ as in~\cite{KW}. 
For a closed subvariety $i: Z\hookrightarrow X$ of dimension $d$ and a smooth open dense subset $j: U\hookrightarrow Z$, we denote by
\[ \delta_Z \;=\; \IC_Z[d] \;=\; i_*j_{!*}(k_U[d]) \in \Perv(X) \] 
the perverse intersection cohomology sheaf of $Z$ and put $\lambda_Z = i_* j_*(k_U[d])$. %
If~$Z$ is normal, $\lambda_Z=k_Z[d]$ by~\cite[lemma 1, p.~5]{WeBN}.
Both~$\delta_Z$ and~$\lambda_Z$ are of geometric origin~\cite[sect.~6.2.4]{BBD} and only depend on $Z$ but not on the choice of $U$. Thus for smooth~$Z$ we get $\delta_Z=\lambda_Z=k_Z[d]$.

\subsection{A semisimple quotient category}\label{subsec:multiplier}
In~\cite[cor.~6, p.~36]{WeBN} it has been shown that every simple perverse sheaf $\delta \in \Perv(X)$ of geometric origin with~$\bbH^\bullet(X,\delta)\neq 0$ is a rigid object in $(\Dbc(X, k), *)$ with dual 
\[ \delta^\vee \;=\; (-\id_X)^*(D(\delta)), \]
where $D$ denotes the Verdier dual.
Using this we want to construct as in~\cite{WeT} a semisimple $k$-linear rigid abelian tensor category $\Pbar(X)$ whose objects are given by semisimple perverse sheaves.
Note that in general the convolution of two perverse sheaves does not have to be perverse. 

\medskip

By~\cite{KrW2} the perverse sheaves $\delta\in \Perv(X)$ with Euler characteristic zero define a thick subcategory $T(X)\subset \Dbc(X, k)$, and a perverse sheaf lies in $T(X)$ iff all its constituents are in $T(X)$. 
Clearly $T(X)$ defines a tensor ideal in $\Dbc(X, k)$ with respect to the convolution product, which implies that the quotient category $\Dbcbar(X, k)=\Dbc(X, k)/T(X)$ is again a $k$-linear tensor category. All simple objects in this latter quotient category are rigid.

\medskip

Although the full abelian subcategory $P(X)\subset \Dbc(X, k)$ of semisimple perverse sheaves on~$X$ is not stable under the convolution product, it is shown in loc.~cit. that its image $\Pbar(X)$ in $\Dbcbar(X, k)$ indeed {\em is} stable under this product, which via \cite{WeT} amounts to saying that every (semisimple) perverse sheaf on $X$ is a multiplier. So $\Pbar(X)$ is a $k$-linear semisimple rigid abelian tensor category under convolution. Similarly, via ~\cite[lemma 10, p.~36]{WeBN} it follows that the mixed perverse sheaves define a $k$-linear rigid abelian tensor category $\Pbar_{mixed}(X)$ under convolution, and this category contains $\Pbar(X)$ as a full abelian tensor subcategory.

\subsection{The symmetric square of the theta divisor} \label{subsec:theta_square} 
Let us now assume that the theta divisor $\Theta \subset X$ is normal and hence irreducible \cite{EL}.
The convolution square $\delta_{\Theta}\!*\!\delta_{\Theta}$ contains the unit object $\one$ precisely once because $\delta_\Theta$ is a simple self-dual object of the rigid symmetric monoidal abelian quotient category $\Pbar(X)$. Furthermore the commutativity constraint 
$S: \delta_\Theta *\delta_\Theta \cong \delta_\Theta*\delta_\Theta$
of~\cite[sect.~2.1]{WeBN} acts as multiplication by $(-1)^{g-1}$ on $\calH^0(\delta_\Theta*\delta_\Theta)_0$. Indeed $\delta_\Theta = \IC_\Theta[g-1]$, and the commutativity constraint for the intersection cohomology sheaf $\IC_\Theta$ is the identity on $\calH^{2g-2}(\IC_\Theta*\IC_\Theta)_0$. This follows from considering fundamental classes, since we may replace $\IC_\Theta$ by the constant sheaf as $\Theta$ is normal.  Notice, the shift by $g-1$ accounts for the factor $(-1)^{g-1}$.

\medskip

This being said, it follows that $\one$ lies in the alternating square~$\Lambda^2(\delta_\Theta)$ for even~$g$ and 
in the symmetric square $S^2(\delta_\Theta)$ for odd $g$.
So there are perverse sheaves $\delta_\pm$ without constituents from $T(X)$ and complexes $\tau_\pm \in T(X)$ such that
\[
 S^2(\delta_{\Theta}) \;=\;
 \left\{
 \begin{array}{c}
 \delta_+ \oplus \tau_+ \\
 \one\oplus \delta_+ \oplus \tau_+
 \end{array}
 \right.
 \textnormal{and} 
 \quad
 \Lambda^2(\delta_{\Theta}) \;=\; 
 \left\{
 \begin{array}{c} \one \oplus \delta_- \oplus \tau_- \\
 \delta_- \oplus \tau_-
 \end{array}
 \right.
 \textnormal{for}
 \;
 \left\{
\begin{array}{c} \textnormal{$g$ even},\\
 \textnormal{$g$ odd}.
 \end{array}
 \right.
\]

If $\Theta$ is smooth, our construction of the family $\pi: Y_U\to U$ in section~\ref{subsec:family} shows that 
$R\pi_*(k_Y)[2g-2]$ corresponds to $(\delta_\Theta * \delta_\Theta)|_U$.
To prove our claim from the introduction that conjecture~\ref{conj:thetasquare} implies conjecture~\ref{conj:vhs}, let us check that $\calV_{\pm}[2g-2]$, considered as Hodge modules in the sense of \cite{SaMHP} and \cite{SaIntro},
have the underlying perverse sheaves $\delta_{\pm}|_U$.

\medskip

Indeed, the commutativity constraint is $S = Ra_*(\phi)$
for the involution~$\phi$ of $k_\Theta[g-1]\boxtimes k_\Theta[g-1]$ given by~$\phi(s\boxtimes t) = (-1)^{g-1} \cdot t\boxtimes s$.
Thus $S|_U$ is the~$\sigma$ of section~\ref{subsec:family} on $R\pi_*(k_{Y_U})$ twisted by~$(-1)^{g-1}$, and  $S^2(\delta_\Theta)|_U$ and~$\Lambda^2(\delta_\Theta)|_U$ are the part of $R\pi_*(k_{Y_U})$ on which $\sigma$ acts by~$\pm (-1)^{g-1}$ respectively.
It only remains to note that the above complexes $\tau_\pm$ are the constant subvariations of section~\ref{subsec:vhs}, which is a consequence of the following lemma.

\subsection{Lemma}  \label{lem:tau}
{\em If $\Theta$ is smooth, then
\[
 \tau_{\pm}  \; = \; \bigoplus_{\textnormal{\em $\mu$ odd for ``$+$''}\atop\textnormal{\em $\mu$ even for ``$-$''}}^{} H^{g-2-|\mu|}(X, k) \otimes \delta_X[\mu].
\]
}

\medskip

{\em Proof.} By semisimplicity we have a decomposition $\tau_\pm = \tau_\pm' \oplus \tau_\pm''$ where~$\tau_\pm'$ denote the direct sum of all complex shifts of $\delta_X$ that enter $\tau_\pm'$. In particular we then have $\bbH^\bullet(X, \tau_\pm'')=0$ since every translation-invariant simple perverse sheaf different from $\delta_X$ has vanishing hypercohomology~\cite[sect.2.3]{WeBN}. Hence $\tau_\pm'$ can be computed from hypercohomology as follows:

\medskip

Using the K\"unneth formula $\bbH^\bullet(X, \delta_\Theta*\delta_\Theta)=\bbH^\bullet(X, \delta_\Theta)\otimes \bbH^\bullet(X, \delta_\Theta)$ and the fact that $\bbH^\nu(X, \delta_\pm)=0$ for $\nu \leq -g$ by lemma~\ref{subsec:lowest-hypercohomology}, one then checks that in perverse cohomology degrees $\mu \leq 0$ the complexes $\tau_\pm'$ coincide with the right hand side of the lemma. By the hard Lefschetz theorem the same then also holds in perverse cohomology degrees $\mu>0$. Finally, using the result for $\tau_\pm'$ for $\mu=0$, we have the stalk cohomology
\[
 \calH^{-g}(\tau_\pm')_0 = H^{g-2}(X, k) = H^{g-2}(\Theta, k) = \calH^{-g}(\delta_\Theta*\delta_\Theta)_0 = \calH^{-g}(\tau_\pm'\oplus \tau_\pm'')_0\ ,
\]
so the non-constant translation-invariant complexes $\tau_\pm''$ must be zero.
\qed

\section{Tannakian categories (conjecture~\ref{conj:thetagroup})}  \label{sec:super-tannakian}

We now discuss the construction the algebraic group $G(X,\Theta)$ that occurs in conjecture~\ref{conj:thetagroup}, again following~\cite{WeT} and \cite{KrW2}.

\subsection{The Tannaka group attached to the theta divisor}\label{subsec:G_X_Theta} 
Suppose now that the theta divisor $\Theta \subset X$ is normal. Inside the rigid abelian tensor category $\Pbar(X)$ of section~\ref{subsec:multiplier}, 
we consider the full abelian tensor subcategory $ \Pbar(X, \Theta) \;=\; \langle \delta_\Theta \rangle$ generated by the perverse sheaf $\delta_\Theta$.
By the general Tannkian description in~\cite{KrW2} there exists an affine algebraic group $G=G(X, \Theta)$ over $k$ and an equivalence
\[ \omega: \; \Pbar(X, \Theta) \; \stackrel{\cong}{\longrightarrow} \; Rep(G) \]
of tensor categories, where $Rep(G)$ denotes the category of algebraic representations of $G$ over $k$. Since $\delta_\Theta$ is a tensor generator of $\Pbar(X,\Theta)$, the action of $G(X,\Theta)$ on~$\omega(\delta_\Theta)$ is faithful. Let $P(X,\Theta) \subset D_c^b(X,k)$ denote the full monoidal subcategory which is the preimage of the quotient category $\Pbar(X,\Theta)$.

\subsection{Euler characteristics}  \label{subsec:sdim}
Since in our quotient category we only divided out perverse sheaves of Euler characteristic zero, for objects $\gamma \in \Pbar(X, \Theta)$ the Euler characteristic $\chi(\gamma)=\sum_{i\in \bbZ} (-1)^i \dim_k(\bbH^i(X,\gamma))$ is well-defined. We claim that
\[ \dim_k (\omega(\gamma)) \; = \; \chi(\gamma). \]
Indeed, the composite
$\varphi: \; \one \stackrel{coev}{\longrightarrow} \gamma*\gamma^\vee \cong \gamma^\vee*\gamma \stackrel{ev}{\longrightarrow} \one$ of coevaluation and evaluation in $\Pbar(X, \Theta)$ satisfies
 $\omega(\varphi) = d\cdot \id$ for $d = \dim_k (\omega(\gamma))$. Hence
$\bbH^\bullet(X, \varphi) = d\cdot \id$, and since $\bbH^\bullet(\varphi)$ is also the composite of coevaluation and evaluation in the category of super vector spaces, we get $d = \chi(\gamma)$.

\subsection{An upper bound} \label{subsec:thetagroupmotivation}
On the Tannakian side, having the unit object $\one$ in $\Lambda^2(\delta_\Theta)$ resp.~$S^2(\delta_\Theta)$ means that the representation $\omega(\delta_\Theta)$ preserves some symplectic resp. orthogonal bilinear form. So we obtain

\begin{lem*}
If $\Theta$ is smooth, then 
\[
 G(X, \Theta) \; \subseteq \;
 \begin{cases}
  \Sp(g!, k) \; \; \textnormal{\em for $g$ even}, \\
  \SO(g!, k) \; \; \textnormal{\em for $g$ odd}.
 \end{cases}
\]
\end{lem*}

{\em Proof.} By the dimension formula in section~\ref{subsec:sdim}, a Gauss-Bonnet calculation shows
$\dim_k (\omega(\delta_\Theta))= g!$ if $\Theta$ is smooth.
Furthermore, in section~\ref{subsec:theta_square} we have seen that $\one$ occurs in $S^2(\delta_\Theta)$ resp.~$\Lambda^2(\delta_\Theta)$ for $g$ odd resp.~even. This proves the lemma with $\mathrm{O}(g!, \bbC)$ in place of $\SO(g!, \bbC)$. But for a generic ppav it follows from~\cite{KrW2} that the group $G(X, \Theta)$ does not admit non-trivial characters, hence it must be contained in $\SO(g!, k)$ if $g$ is odd. To deduce the result for arbitrary ppav's with a smooth theta divisor one can then use a specialization argument. \qed
\\

This lemma and the cases $g=2,3,4$ motivate conjecture~\ref{conj:thetagroup}.
For~$g=2$ every ppav $(X, \Theta)$ with smooth $\Theta$ is the Jacobian of a hyperelliptic curve; this case is covered by~\cite[p.~124 and th.~14]{WeBN}.
For $g=3$ every ppav $(X, \Theta)$ with smooth $\Theta$ is the Jacobian of a non-hyperelliptic curve, and $\omega(\delta_\Theta)$ is the second fundamental representation of $G(X, \Theta) = \Sl(4, \bbC)/\mu_2$ by loc.~cit. It corresponds to the standard representation of $\mathrm{SO}(6,\bbC)$ via the exceptional isomorphism $\Sl(4,\bbC)/\mu_2 \cong \mathrm{SO}(6, \bbC)$.

\section{Tensor functors defined by nearby cycles} \label{sec:cycles}

To study the tensor categories $\Pbar(X)$ or $\Pbar(X, \Theta)$ of section~\ref{subsec:multiplier} when~$(X, \Theta)$ varies in families, we briefly recall some facts from nearby cycle theory as exposed in~\cite[exp.~XIII-XIV]{SGA7} \cite{KS}.
Let
$f: \calX \rightarrow S$ be a proper surjective algebraic morphism from a smooth complex algebraic variety to a smooth algebraic curve. For $s\in S(\bbC)$ we want to relate complexes on the special fibre $\calX_s=f^{-1}(s)$ to complexes on $\calX_t=f^{-1}(t)$ for~$t$ in some pointed analytic neighborhood of~$s$.

\subsection{Analytic nearby cycles}\label{subsec:nearby} 
Let $D\subset S(\bbC)$ be a small coordinate disc centered at $s$. The morphism $f: \calX \to S$ yields a proper holomorphic map
$\calX_D = \calX \times_S D \to D$ of analytic spaces,
ditto with $D$ replaced by the pointed disc $D^*=D\setminus \{s\}$ or by the universal covering $\tilde{D}^*$ of $D^*$. So we have a cartesian diagram
\[
 \xymatrix@M=0.5em{
 \calX_s \ar@{^{(}->}[r]^i \ar[d] & \calX_D \ar[d] & \calX_{D^*} \ar@{_{(}->}[l]_j \ar[d] & \calX_{\tilde{D}^*} \ar[l]_\pi \ar[d]  \\
 \{s\} \ar@{^{(}->}[r] & D & D^* \ar@{_{(}->}[l] & \tilde{D}^* \ar[l]
 }
\]
where $\pi$ is a covering map and $i$ resp.~$j$ are closed resp.~open immersions.
Let $\tilde{j}=j\circ \pi$. Following~\cite[exp.~XIV]{SGA7} and~\cite[sect.~8.6]{KS} we consider the functor
\[
 \Psi: \; \Dbc(\calX, \bbC) \longrightarrow \Dbc(\calX_s, \bbC), \;
 \delta\mapsto i^* \, R\tilde{j}_* \, \tilde{j}^* \, (\delta|_{\calX_D})
\]
of nearby cycles. It factors over a functor $\Dbc(\calX_{D^*}, \bbC) \to \Dbc(\calX_s, \bbC)$ which we also denote by $\Psi$. If $D$ has been chosen sufficiently small, one has an isomorphism
\begin{equation} \label{aha}
\bbH^i(\calX_s, \Psi(\delta)) \; \cong \; \bbH^{i}(\calX_t, \delta|_{\calX_t})
\quad
\textnormal{for all}
\quad
t\in D^*.
\end{equation}
For $\delta\in \Dbc(\calX_D, \bbC)$ the pullback of the morphism $\delta\rightarrow \tilde{j}_*\tilde{j}^*(\delta)$ under $i$ defines a morphism $i^*(\delta)\rightarrow \Psi(\delta)$. As in~\cite[eq.~8.6.7]{KS} the cone of this morphism defines the functor $
 \Phi:\; \Dbc(\calX_D, \bbC) \rightarrow \Dbc(\calX_s, \bbC)
$
of vanishing cycles, so for every $\delta\in \Dbc(\calX_D, \bbC)$ we have a distinguished triangle\label{subsec:psiperv}
\begin{equation} \label{iPsiPhiTriangle}
 i^*(\delta) \longrightarrow \Psi(\delta) \longrightarrow \Phi(\delta) \longrightarrow i^*(\delta)[1] 
\end{equation} 
in $\Dbc(\calX_s, \bbC)$. Note the various shifting conventions in the literature;~our $\Phi$ would be denoted $\Phi[-1]$ in~\cite{KS}. With our conventions, the shifted functors $\Psi[-1]$ and $\Phi[-1]$ commute with Verdier duality by~\cite[1.4]{Br}, and they send $\Perv(\calX)$ to $\Perv(\calX_s)$ in view of~\cite[cor.~10.3.13]{KS}.

\subsection{Algebraic nearby cycles} \label{subsec:dvr}
Localizing in the point $s\in S(\bbC)$, let us now replace $S$ by the spectrum of a Henselian discrete valuation ring centered at
the special point $s$, and denote by~$\eta$ its generic point. We then have an algebraic version of the nearby cycles, a functor $\Psi: \Dbc(\calX) \to \Dbc(\calX_s \times_s \eta)$ in the sense of \cite[\S 3.1 and \S 4]{Il}. It factors over a functor $\Dbc(\calX_\eta) \to \Dbc(\calX_s \times_s \eta)$ which we also denote by $\Psi$. The properties described in the analytic setting carry over to this case. We have a sequence (\ref{iPsiPhiTriangle}), and $\Psi$ maps $\Perv(\calX_\eta)$ to $\Perv(\calX_s)$ and commutes with Verdier duality by th.~4.2 and cor.~4.5 of loc.~cit. Furthermore, if $\calX\to S$ is an abelian scheme, then $\Psi$ is a tensor functor with respect to convolution on $\calX_\eta$ resp.~$\calX_s$ by compatibility with proper maps~\cite[exp.~XIII, 2.1.7]{SGA7} and by the K\"unneth formula~\cite[5.1]{BB}.

\subsection{Group-theoretical reformulation} \label{subsec:grouptheoretical}
For an abelian scheme $\calX \to S$, in the setting of~\ref{subsec:dvr}, let $\calT$ be a tensor subcategory of $\Pbar_{mixed}(\calX_\eta)$. Denote by~$\calT_\Psi$ be the tensor subcategory of $\Pbar_{mixed}(\calX_s)$ generated by the image $\Psi(\calT)$. If $\calT$ is a finitely generated tensor category, so is  $\calT_\Psi$. Then there are algebraic $k$-groups $G$ and $G_\Psi$ such that ${\calT} = \Rep(G)$ and $\calT_\Psi = \Rep(G_\Psi)$, where the right hand sides denote the tensor categories of algebraic representations of $G$ resp.~$G_\Psi$.

 \begin{lem*}  The functor $\Psi$ is a $k$-linear $\otimes$-functor ACU and maps perverse sheaves  in $T(\calX_\eta)$ to perverse sheaves in $T(\calX_s)$. Hence it induces a $k$-linear exact $\otimes$-functor
 $$ \Psi: \Rep(G) \longrightarrow \Rep(G_\Psi).$$ 
 \end{lem*} 

{\it Proof}. 
The algebraic analog of the isomorphism (\ref{aha}) shows that $\Psi$
maps complexes with vanishing Euler characteristic to  complexes with vanishing Euler characteristic. So we get a $\otimes$-functor $$\Pbar_{mixed}(\calX_\eta)\longrightarrow \Pbar_{mixed}(\calX_s),$$ which immediately implies the assertions. Notice that $\Psi$ maps distinguished triangles to distinguished triangles, hence induces an exact functor. \qed 

\medskip

Recall that Deligne \cite[sect.~8]{DelTa} has attached to any Tannaka category~$\calT$ an $\mathrm{Ind}(\calT)$-groupscheme $\pi(\calT)$, called the fundamental group of~$\calT$. By 8.15 of loc.~cit. any $k$-linear exact $\otimes$-functor $\eta: \calT_1 \to \calT_2$
induces a morphism
\begin{equation}\label{Deligne}
\pi(\calT_2) \;\longrightarrow \; \eta(\pi(\calT_1)).
\end{equation}
Under the weak conditions (2.2.1) and (8.1) of loc.~cit. (which are verified for representation categories $\calT_i =\Rep(G_i)$ of algebraic groups~$G_i$ over an algebraically closed field $k$ of characteristic zero) theorem 8.17 of loc.~cit. implies that the functor $\eta$ induces an equivalence of $\calT_1$ with the category of
objects in $\calT_2$ endowed with an action of $\eta(\pi(\calT_1))$ such that the natural action of $\pi(\calT_2)$ is induced by (\ref{Deligne}) above. If $\eta$ is a fiber functor
to the tensor category $\calT_2=Vec_k$ of finite-dimensional vector spaces over $k$, this reduces to the assertion $\calT_1 = Rep(G_1)$ for $G_1=\eta(\pi(\calT_1))$.

\subsection{Subgroups defined by degenerations} \label{subsec:groupinclusion}
Let $\calT_1 = \calT$ be a finitely generated tensor subcategory of $\Pbar_{mixed}(\calX_\eta)$ as above, and choose a fiber functor $\omega$ of $\calT_2=\calT_\Psi$. Then $\omega\circ\Psi$ is a fiber functor of $\calT_1=\calT$, since it is exact and therefore faithful by the isomorphism~(\ref{aha}). Hence~(\ref{Deligne}) applied to the functor $\Psi$ induces a morphism of algebraic $k$-groups
 $$ G_\Psi = G_2 \; \longrightarrow \; G_1 = G.$$
In \cite[p.~118]{DM} it is shown that this morphism $G_\Psi \to G$ 
is a closed immersion iff every object $K$ of $\calT_2 = Rep(G_\Psi)$ is isomorphic to a subquotient of an object $\Psi(K')$ for some~$K'$ in $\calT_1 = \Rep(G)$. In our situation this holds by the definition of $\calT_\Psi$, so we get the 

\begin{lem*} 
The algebraic $k$-group $G_\Psi$ is a closed algebraic $k$-subgroup of $G$, and $\Psi$ can be identified with the restriction functor $\Rep(G)\to \Rep(G_\Psi)$.
\end{lem*}

\subsection{Semisimplification}
In dealing with Tannakian categories of mixed perverse sheaves on abelian varieties $X$ over $k$ one can use the following 

\begin{lem*}
Let $\langle \delta \rangle = \Rep(G)$ be the full tensor subcategory of~$\Pbar_{mixed}(X)$ which is generated by a mixed perverse sheaf~$\delta$. Then the full tensor subcategory generated by the semisimplification $\delta^{ss}$ is
\[ \langle \delta^{ss} \rangle = \Rep(G^{red}) \]
where $G^{red} = G / R_u(G)$ denotes the quotient of $G$ by its unipotent radical. 
\end{lem*}

{\em Proof.} This boils down to a statement about the categories of representations of algebraic groups over a field of characteristic zero; see~\cite{KrW2}.
\qed

\label{subsec:psi-cohomology}

\section{Local Monodromy} \label{subsec:monodromy} 

In the setting of~\ref{subsec:nearby}, let $\delta\in \Perv(\calX)[-1]$. Then $\Psi(\delta)$ is a perverse sheaf on $\calX_s$, and for fixed $t\in D^*$ the action of $\pi_1=\pi_1\!(D^*, t)$ on the universal cover~$\tilde{D}^*$ induces a monodromy operation on this perverse sheaf. In the algebraic setting of section~\ref{subsec:dvr} we can proceed similarly, replacing the fundamental group $\pi_1$ by the pro-cyclic local monodromy group $\bbZ_l(1)$ as in~\cite[\S 3.6]{Il}, see section~\ref{lem:GrTensor} below.

\subsection{Unipotent nearby cycles}
Let $T$ be a generator of $\pi_1$ acting on $\Psi(\delta)$ as above. 
We have a direct sum decomposition 
\[ \Psi(\delta) \;=\; \Psi_1(\delta) \oplus \Psi_{\neq 1}(\delta) \]
where $\Psi_1(\delta)\subseteq \Psi(\delta)$ denotes as in~\cite[lemma 1.1]{Rei} the maximal perverse subsheaf on which $T$ acts unipotently. Similarly $\Phi(\delta) = \Phi_1(\delta)\oplus \Phi_{\neq 1}(\delta)$. We say that~$\delta$ has unipotent global monodromy if $\bbH^\bullet(\calX_s, \Psi_1(\delta)) = \bbH^\bullet(\calX_s, \Psi(\delta))$. Then from the Picard-Lefschetz formulas \cite[exp.~XV, th.~3.4(iii)]{SGA7} one draws the

\begin{lem*}\label{Lef}
If $\calX_s$ is regular except for finitely many ordinary double points and if $\dim(\calX_s)$ is odd, then $\delta=\delta_\calX[-1]$ has unipotent global monodromy; more precisely $(T-1)^2$ acts trivially on $\bbH^\bullet(\calX_s, \Psi(\delta))$. \label{rem:cone} 
\end{lem*}

\label{subsec:sp_sequence}\label{rem:sp_weight} 
Returning to the general case, since by definition $T-1$ acts nilpotently on~$\Psi_1(\delta)$, we can define
$N = \frac{1}{2\pi i} \log(T):  \Psi_1(\delta) \longrightarrow \Psi_1(\delta)(-1)$. The cone of~$N$ in $\Dbc(\calX_s, \bbC)$ is given by
\[  C\bigl(\Psi_1(\delta) \stackrel{N}{\longrightarrow} \Psi_1(\delta)(-1)\bigr) \; = \;
 C\bigl(\Psi(\delta) \stackrel{\frac{T-1}{2\pi i}}{\longrightarrow} \Psi(\delta)(-1)\bigr) \; = \; i^* Rj_*j^*(\delta[1]). \]
Indeed, the first equality holds because $T-1$ is an isomorphism on $\Psi_{\neq 1}(\delta)$ whereas on $\Psi_1(\delta)$ its kernel and cokernel coincide with those of $N$ up to a weight shift. For the second equality see~\cite[eq.~(3.6.2)]{Il} and the remarks thereafter. The perversity of $\Psi_1(\delta)$ and the above formula for the cone of $N$ imply that if we define specialization functors by
\[
 sp(-) \;=\; \pH^0(i^*Rj_*j^*(-)) \quad \textnormal{and} \quad sp^\dag(-) \;=\; \pH^1(i^*Rj_*j^*(-)),
\]
we obtain an exact sequence of perverse sheaves on $\calX_s$
\[
 0 \longrightarrow sp(\delta) \longrightarrow \Psi_1(\delta) \stackrel{N}{\longrightarrow} \Psi_1(\delta)(-1) \longrightarrow sp^\dag(\delta) \longrightarrow 0.
\]
Since $\Psi$ and hence also $\Psi_1$  preserve distinguished triangles, the functor $sp$ is left exact on perverse sheaves. 

\subsection{The monodromy filtration on $\Psi_1(\delta)$}\label{MF}
As in~\cite[section 1.6]{DelW} the nilpotent operator $N$ gives rise to a unique finite increasing filtration~$F_\bullet$ of~$\Psi_1(\delta)$ in $\Perv(\calX_s)$ such that for all $i$,
\begin{itemize}
 \item $N(F_i(\Psi_1(\delta)))\subset F_{i-2}(\Psi_1(\delta))(-1)$, and
 \item $N^i$ induces an isomorphism $Gr_i(\Psi_1(\delta)) \stackrel{\cong}{\rightarrow} Gr_{-i}(\Psi_1(\delta))(-i)$.
\end{itemize}
Each $Gr_{-i}(\Psi_1(\delta))$ with $i\geq 0$ has an increasing filtration with composition factors $P_{-i}(\delta)$, $P_{-i-2}(\delta)(-1)$, $P_{-i-4}(\delta)(-2)$, $\dots$ where
\[ P_i(\delta) := \ker(N: Gr_i(\Psi_1(\delta))\rightarrow Gr_{i-2}(\Psi_1(\delta))(-1)), \]
In what follows we will represent this situation as in loc.~cit.~by a triangle
\[
 \xymatrix@=0.1em{
 \vdots & \quad & & & & \quad \dots \\
 Gr_2(\Psi_1(\delta)) & \quad & & & { P_{-2}(\delta)(-2)} \ar[dd]^N_\cong &  \\
 Gr_1(\Psi_1(\delta)) & \quad & & { P_{-1}(\delta)(-1)} \ar[dd]^N_\cong & & \quad \dots \\
 Gr_0(\Psi_1(\delta)) & \quad & { \quad P_0(\delta) \quad } & & { P_{-2}(\delta)(-1)} \ar[dd]^N_\cong & \\
 Gr_{-1}(\Psi_1(\delta)) & \quad & & { \quad P_{-1}(\delta) \quad } & & \quad \dots \\
 Gr_{-2}(\Psi_1(\delta)) & \quad & & & { \quad P_{-2}(\delta) \quad } & \\
 \vdots & \quad & & & & \quad \dots
 }
\]
where each line gives the decomposition of the corresponding graded piece. 
The lower boundary entries $P_0(\delta), P_ {-1}(\delta), P_{-2}(\delta), \dots$ in the triangle are the graded pieces of $sp(\delta) = \ker(N)$, with $P_0(\delta)$ as the top quotient. In the situation of proposition~\ref{cor:sp}(a) below, the entries above these lower boundary entries belong to~$\Phi_1(\delta)$.

\subsection{The monodromy filtration on $\Psi(\delta)$}
For any $\delta \in \Perv(\calX)[-1]$, the local monodromy theorem~\cite[th.~2.1.2]{Il} and the Jordan decomposition of the nearby cycles~\cite[lemma~4.2]{Rei} show that there is an $a\in \bbN$ such that $T^a - 1$ is nilpotent on all of $\Psi(\delta)$. In this case, using the operator $N' := \frac{1}{2\pi i} \log(T^a): \Psi(\delta)\to \Psi(\delta)(-1)$ in place of $N: \Psi_1(\delta)\to \Psi_1(\delta)(-1)$ one can define a filtration $F'_\bullet$ as in section~\ref{MF} on all of $\Psi(\delta)$. This filtration does not depend on the chosen~$a$ such that $T^a-1$ is nilpotent. Even though in general one has $N'|_{\Psi_1(\delta)}\neq N$, the fact that $T$ acts unipotently on $\Psi_1(\delta)$ implies that the kernel and the image of $N'|_{\Psi_1(\delta)}$ are the same as those of $N$. Hence 
\[
 F_\bullet(\Psi_1(\delta)) \;=\; \Psi_1(\delta) \cap F_\bullet'(\Psi(\delta)).
\]
Notice however that $sp(\delta) = \ker(N: \Psi_1(\delta)\to \Psi_1(\delta)(-1))$, as defined in~\ref{subsec:sp_sequence}, will in general only be a perverse {\em subsheaf} of $\ker(N': \Psi(\delta)\to \Psi(\delta)(-1))$ because $N'$ may have a non-trivial kernel on~$\Psi_{\neq 1}(\delta)$. On the other hand, working with $\Psi(\delta)$ instead of $\Psi_1(\delta)$ has the following advantage.

\subsection{Tensor functoriality} 
All of the above has an analog in the algebraic setting of section~\ref{subsec:dvr}, if $T$ is a topological generator of the local monodromy group $\bbZ_l(1)$ as in~\cite[\S 3.6]{Il}.

\begin{lem*} \label{lem:GrTensor}
For $\delta \in \Perv(\calX_\eta)$, denote by $Gr_\bullet'(\Psi(\delta)) = \bigoplus_{i\in \bbZ} Gr_i'(\Psi(\delta))$ the associated graded with respect to the filtration $F'_\bullet$ on $\Psi(\delta)$ as defined above. Then we have an induced functor
\[
 Gr_\bullet'\circ \Psi: \; \Pbar(\calX_\eta) \; \longrightarrow \; \Pbar(\calX_s), \; \delta \mapsto Gr_\bullet'(\Psi(\delta))
\]
which is a tensor functor with respect to convolution.
\end{lem*}

{\em Proof.} By lemma~\ref{subsec:grouptheoretical}, $\Psi$ induces a tensor functor $\Pbar(\calX_\eta)\to \Pbar_{mixed}(\calX_s)$ with respect to convolution. If a fixed generator of $\bbZ_l(1)$ acts on $\delta_i\in \Perv(\calX_\eta)$ via endomorphisms $T_i: \Psi(\delta_i)\to \Psi(\delta_i)$ for $i\in \{1,2\}$, then this generator acts on the convolution product $\Psi(\delta_1*\delta_2)=\Psi(\delta_1)*\Psi(\delta_2)$ via the product $T_1*T_2$. Since the tensor subcategory of $\Pbar_{mixed}(\calX_s)$ generated by $\Psi(\delta_1)$ and $\Psi(\delta_2)$ is equivalent to the category of representations of some algebraic group, as in \cite[prop.~1.6.9]{DelW} one deduces $Gr_i'(\Psi(\delta_1)*\Psi(\delta_2)) = \bigoplus_{i_1+i_2=i} Gr_{i_1}'(\Psi(\delta_1))*Gr_{i_2}'(\Psi(\delta_2))$.
\qed

\subsection{Weights}
Concerning weights we have the following observation which is due to Gabber~\cite[1.19]{SaIntro} \cite[th~5.1.2]{BB}. 

\begin{lem*}
If $\delta$ underlies a pure Hodge module of weight~$w$ on $\calX_\eta$, then each graded piece $Gr_i'(\Psi(\delta))$ is pure of weight $w+i$ and the filtration $F_\bullet'(\Psi(\delta))$ is the weight filtration of $\Psi(\delta)$ up to an index shift. Hence the analogous statement also holds for the filtration $F_\bullet$ of $\Psi_1(\delta)$ .
\end{lem*}

\section{Behavior of $G(X, \Theta)$ under specialization}
\label{7}

In this section we study the behavior of the Tannaka group $G(X, \Theta)$ when the ppav $(X, \Theta)$ degenerates. For this we need to control the specialization functor $sp$ as defined in section~\ref{subsec:sp_sequence} above.

\subsection{Proposition.}\ \label{cor:sp}
{\em With notations as in section~\ref{subsec:nearby}, for $\delta\in \Perv(\calX)[-1]$ the specialization is given by
\[ sp(\delta)=i_* (j_{!*} \gamma)[-1]
\quad \textnormal{\em with the perverse sheaf} \quad
\gamma = \delta[1]|_{\calX_{D^*}}.
\]  
This in particular leads to the following two observations.

\begin{enumerate}
 \item If $\delta=(j_{!*}\gamma)[-1]$, then $sp(\delta)=i^*(\delta)$. In this case the triangle~(\ref{iPsiPhiTriangle}) yields an exact sequence of perverse sheaves
\[
 0 \longrightarrow sp(\delta) \longrightarrow \Psi(\delta) \longrightarrow \Phi(\delta) \longrightarrow 0,
\]
and the same also holds with $\Psi$ and~$\Phi$ replaced by $\Psi_1$ and $\Phi_1$. 

\item If $\delta$ underlies a mixed Hodge module of weight $\leq w$, so does~$sp(\delta)$. 
\end{enumerate}
}

\medskip

{\em Proof.} Part {\em (a)} follows directly from the formula for $sp(\delta)$. For exactness of the sequence it does not matter whether we use the functors $\Psi$ and $\Phi$ or their unipotent counterparts because the monodromy operator $T$ acts trivially on $sp(\delta)$. Part {\em (a)} also follows directly from the formula for the specialization, together with the permanence properties of weights under pull-back and intermediate extensions. It remains to show that
\[ sp(\delta)  \; = \;   i^* (j_{!*} \gamma )[-1].
 \]
For this notice that a basic property of intermediate extensions \cite[III.5.1]{KW} implies $i^*j_{!*}\gamma = \ptau_{<0}i^*Rj_*\gamma$. Furthermore we have $\ptau_{<0}i^*Rj_*\gamma = sp(\delta)[1]$ since $j$ is affine so that $i^*Rj_*\gamma \in \pD^{[-1,0]}(\calX_0)$ in view of Artin's vanishing theorem~\cite[th.~4.1.1]{BBD}. See also the formula for the cone of $N$ in the proof of lemma~\ref{subsec:sp_sequence}. \qed

\subsection{Some further remarks}
In the situation of part (a) above, the local invariant cycle theorem~\cite[cor.~6.2.9]{BBD} states that for all $i\in \bbZ$ the  induced morphisms
$\bbH^i(\calX_s, sp(\delta)) \to \bbH^i(\calX_s, \Psi(\delta))$ factor through epimorphisms
\[
 \bbH^i(\calX_s, sp(\delta)) \; \twoheadrightarrow \; \bbH^i(\calX_s, \Psi(\delta))^T
\]
onto the invariants under the local monodromy group. In (b), if $\delta$ is pure of some weight $w$, we denote by $\spbar(\delta)$ the highest top quotient (of weight $w$) of the weight filtration of $sp(\delta)$.

\begin{ex*} \label{ex:sp}
(i) If $\calX$ is smooth over $S$ of relative dimension $d$, then
$$ sp(\bbC_\calX[d]) = \bbC_{\calX_s}[d]. $$
(ii) If $\calX$ is regular of dimension $d+1$, then $sp(\bbC_{\calX}[d]) = \bbC_{\calX_s}[d]$. If in addition $\calX_s$ is normal and if the only singularities of $\calX_s$ are finitely many ordinary double points, then $\spbar(\bbC_{\calX}[d]) = \delta_{\calX_s}$.
\end{ex*}

{\em Proof.} Part (i) is obvious. To check the first statement in (ii), note that the regularity of $\calX$ implies $j_{!*}j^*(\bbC_{\calX}[d+1]) = \bbC_{\calX}[d+1]$. For the second statement in (ii) see lemma~\ref{subsec:doublepoints}, noting that in the case at hand $\calX_s$ is automatically a local complete intersection. \qed

\medskip

For a principally polarized abelian scheme $\calX \to S$ with relative theta divisor $\Theta_\calX \subset \calX$ over a discrete valuation ring $S$ we denote by $\Theta_s=\Theta_{\calX_s}$ and $\Theta_\eta=\Theta_{\calX_\eta}$ the theta divisors of the special resp.~generic fibre. Then we have the following

\begin{lem*} 
If $\Theta_{\calX}$ and $\Theta_s$ are normal, then $\overline{sp}(\delta_{\Theta_\eta})$ contains $\delta_s$. In particular, then $G_s=G(\calX_s,\Theta_s)$ is a subquotient
of $G =G(\calX_\eta,\Theta_\eta)$, and the defining representations of these groups satisfy 
$\dim_k ( \omega(\delta_{\Theta_s}) ) \leq \dim_k ( \omega(\delta_{\Theta_\eta}) )$.
\end{lem*}

{\em Proof.} a) If $\delta_{\Theta_s}$ is a constituent of $\spbar(\delta_{\Theta_\eta})$, then $G_s$  is a quotient of the Tannaka group~$G_\Psi$, and lemma~\ref{subsec:groupinclusion} implies $G_\Psi \hookrightarrow G$. So the second statement of the lemma follows from the first one. It remains to show that, if $\Theta_{\calX}$ and $\Theta_{s}$ are normal, $\delta_{\Theta_s}$ is a constituent of~$\spbar(\delta_{\Theta_\eta})$.

\medskip
b) If $\Theta_\calX$ is normal, lemma~\ref{subsec:doublepoints}(ii) gives an exact sequence of perverse sheaves
$0 \rightarrow \psi_{\Theta_\calX} \rightarrow k_{\Theta_\calX}[g] \rightarrow \delta_{\Theta_\calX} \rightarrow 0$
where $\psi_{\Theta_\calX}$ has weights $<g$. Recall that $i^*\delta_{\Theta_\calX}[-1] = sp(\delta_{\Theta_\eta})$ is perverse by proposition \ref{cor:sp}.
Hence, if we apply $i^*[-1]$ to this sequence, we get an exact sequence of perverse sheaves on~$\calX_s$
$$ 0 \to \pH^0(i^*\psi_{\Theta_\calX}[-1]) \to k_{\Theta_s}[g-1] \to sp(\delta_{\Theta_\eta}) \to \pH^1(i^*\psi_{\Theta_\calX}[-1]) \to 0  \ ,$$
since $\Theta_s$ is a local complete intersection and therefore also $i^*k_{\Theta_\calX}[g-1] = k_{\Theta_s}[g-1]$ is perverse.
Since the restriction functor~$i^*$ preserves upper bounds on weights, the first
term $\pH^0(i^*\psi_{\Theta_\calX}[-1])$ has weights $<g-1$. 

\medskip

c) Lemma~\ref{subsec:doublepoints}(i) applied to the special fiber gives an exact sequence of perverse sheaves
$$ 0 \to \psi_{\Theta_s} \to k_{\Theta_s}[g-1] \to \delta_{\Theta_s} \to 0 $$
where $\psi_{\Theta_s}$ has weights $<g-1$ with pure quotient $\delta_{\Theta_s}$ of weight $g-1$. Thus the perverse sheaf $k_{\Theta_s}[g-1] $ admits $\delta_{\Theta_s}$ as the highest weight quotient. For weight reasons and by the exact sequence in b), this epimorphism factors over the quotient perverse sheaf $sp(\delta_{\Theta_\eta})$, because $\pH^0(i^* \psi_{\Theta_\calX}[-1])$ has weights $<g-1$. Again for weight reasons the epimorphism then also factors over $\spbar(\delta_{\Theta_\eta})$, which means that we get an epimorphism $\spbar(\delta_{\Theta_\eta})\twoheadrightarrow \delta_{\Theta_s}$. \qed

\medskip
In passing we remark that if in the lemma the assumption 
that $\Theta_s$ is normal is replaced  by the assumption that $\Theta_\eta$ is smooth, then
the argument in b) does imply that $k_{\Theta_s}[g-1]$ is a perverse subsheaf of
$sp(\delta_{\Theta_\eta})$. Indeed, then $\psi_{\Theta_\calX}=i_*(\alpha)$ for some perverse sheaf $\alpha$ supported in the  fiber $\Theta_s$, so that $\pH^0(i^*\psi_{\Theta_\calX}[-1])$ vanishes since $i^*\psi_{\Theta_\calX} = \alpha$ is perverse.

\section{Proof of the main theorem for $g=4$: Outline} \label{sec:outline} 
To prove the main theorem from the introduction, let $X$ be a general ppav of genus 4 with a smooth theta divisor $\Theta \subset X$. We want to show that 
up to skyscraper sheaves the symmetric and alternating squares of $\delta_\Theta$ define irreducible representations of $G$, i.e.~that the perverse sheaves $\delta_\pm$ are simple objects in the representation
category $\Rep(G(X,\Theta))=\Pbar(X,\Theta)$.
For this we compare two tensor functors, 

\begin{enumerate}
 \item the global motivic functor $ \MT:  \MHM(X, \Theta) \to \sRep(\MT(X,\Theta))$
 related to the Mumford-Tate group
of~$\Theta$, and \medskip
 \item a restriction functor induced from an embedding $G_\Psi \hookrightarrow G=G(X,\Theta)$ of algebraic groups via the theory of vanishing cycles for a degeneration of~$X$ into the Jacobian $JC$ of a general curve $C$. \medskip
\end{enumerate}
For part (a) recall the realization functor $D^b(\MHM(X))\to \Dbc(X,k)$ from the bounded derived category of the abelian category of mixed Hodge modules. We define $\MHM(X,\Theta)\subset D^b(\MHM(X))$ to be the preimage $P(X,\Theta)\subset  \Dbc(X,k)$. The direct image under the structure morphism $X\to Spec(\bbC)$ induces the desired tensor
functor $\MT$ to the category of $k$-linear finite-dimensional super representations of the Mumford-Tate group $\MT(X,\Theta)$.%

\medskip
For part (b) we let $(X,\Theta)$ degenerate into the Jacobian~$JC$ of a general curve of genus 4. By the formalism of section \ref{7}, the subgroup $G_\Psi \subset G=G(X,\Theta)$ defined by this degeneration contains the group $G(JC,\Theta_{JC})=\Sl(6, \bbC)/\mu_3$ from~\cite{WeBN}. As explained in section \ref{subsec:psitheta}, the representations of~$G$ associated to the perverse sheaf~$\delta_\Theta$ and its tensor square
become reducible under restriction to the subgroup $G_{\Psi}$, but only with very few simple constituents.  This is the first step of the proof.
To proceed further one considers the same curve degeneration from a Hodge-theoretic point of view. 
Notice that the functors in (a) and (b) 
are not directly related to each other, but only indirectly 
via the following diagram, where the vertical tensor functor is the
composition of the realization functor $\MHM(X,\Theta)\to P(X,\Theta)$ with the quotient functor $P(X,\Theta) \to \Pbar(X,\Theta)=P(X,\Theta)/T(X)$:
$$ \xymatrix{ \MHM(X,\Theta) \ar[d] \ar[r]^-{\MT} &    \sRep(\MT(X,\Theta)) \ar[rr] & & \sRep(\MT(JC,\Theta_{JC}))  \cr
\Pbar(X,\Theta) \ar[r] & \Rep(G) \ar[r] & \Rep(G_\Psi) \ar[r]  &  \Rep(G(JC,\Theta_{JC})) \cr } $$
 Unlike the functor $\MHM(X,\Theta) \to \Rep(G(JC,\Theta_{JC}))$, the Hodge-theoretic 
 functor $\MHM(X,\Theta)\to   \sRep(\MT(X,\Theta))$  is non-trivial on $T(X)$. To compare the two functors we have to keep track carefully  of all constituents from $T(X)$ which have been computed in lemma~\ref{lem:tau}. 

\medskip
The key
lemma~\ref{prop:big} is that the two summands
$\delta_\pm$ of the tensor square of $\delta_\Theta$ are simple up to skyscraper sheaves.
To show this we will construct two large simple subobjects $\gamma_\pm \subseteq \delta_\pm$
and show that they coincide with $\delta_{\pm}$ up to skyscraper sheaves.  For the key lemma we compare the decomposition of the objects $\delta_\pm \in \MHM(X,\Theta)$ in the two representation categories $\sRep(\MT(\Theta))$ and $\Rep(G(JC,\Theta_{JC}))$.
The image of~$\delta_\pm$ in $\sRep(\MT(X,\Theta))$ has many irreducible summands. Even after ignoring the terms from $T(X)$ there appear 28
irreducible representations of the Mumford-Tate group in the representation associated to the tensor square of $H^\bullet(X,\delta_\Theta)$;  the complete list can be read off from table~\ref{table:genericHypercoh} in section~\ref{app:hypercoh}.   

\medskip
For a comparison with (b) the behavior of the Mumford-Tate group under the degeneration of  $(X,\Theta)$ 
into $(JC,\Theta_{JC})$ is relevant. Since $\MT(\Theta_{JC})_{sc} = \Sp(8,\bbC)$ for  a general curve $C$ of genus~4, this behavior is encoded in a group homomorphism
$$ \xymatrix{ \MT(JC,\Theta_{JC})_{sc} = \Sp(8,\bbC) \ar[rr]^-{(id, \varphi)} & &  \Sp(8,\bbC) \times \Sp(10,\bbC) = \MT(X,\Theta)_{sc}} $$
where the homomorphism $\varphi$ is defined in section~\ref{subsec:Bdegen}. 
It is also important to keep track on the natural weight filtrations induced by the curve degeneration
for both functors in situation (a) and (b) together with the action of the local monodromy operator $N$
and the involution $\sigma$.

\medskip
The punchline of the argument is that the local monodromy action 
is nontrivial only for 10 of the 28 cohomology summands, namely those where the representation associated to $B$ is
involved. On each of these 10 summands the local monodromy can be easily described in terms of $B$. In fact, by the local invariant cycle theorem the decomposition related to the functor in (b) is compared to the invariants of the local monodromy on the cohomology, which can be
computed from the invariants $B^N \subset B$ and table~\ref{table:genericHypercoh}.  
If the $N$-invariants $S^N$ of one of the relevant 10 summands $S$ contributes to the decomposition of $\gamma_{\pm}$ induced by the functor in (b), then already the complete summand $S$ is contained in the cohomology of $\gamma_{\pm}$. Since $B^N$ is a large subspace in $B$, this allows
to identify $S$ in the cohomology of the factor $\gamma_{\pm}$ by a comparison with table \ref{table:specialHypercoh}.

\section{Proof of the main theorem for $g=4$: The details}
\label{sec:decomposition}

In order to apply the motivic results of section~\ref{sec:motivic}, we give the proof in an analytic framework. Notice however that all steps in the proof work in the same way if one replaces the Mumford-Tate group $\MT(B)$ by its motivic counterpart $M(B)$ in the sense of Andr\'e using section \ref{motivic}; in particular, the objects $\delta_\pm$ will be shown to be simple not only as mixed Hodge modules but also as perverse sheaves (as required for the main theorem).

\medskip

Recall that if~$V$ is a complex algebraic variety, $\Dbc(V, k)$ is a full subcategory of the derived category of complexes on the analytic space~$V^{an}$ which are constructible for an {\em algebraic} stratification~\cite[6.1.2]{BBD}. Hence it is a full subcategory of the category of all $\bbC$-constructible complexes on~$V^{an}$ in the sense of~\cite{KS}. The number of constituents of a semisimple complex is the same in any of the above triangulated categories in which it lies, so we will denote all of them indifferently by $\Dbc(-, \bbC)$.

\subsection{Degeneration into a Jacobian} \label{subsec:familyofppavs} 
We first construct the degenerating family of ppav's to be used in the proof. Let~$C$ be a generic smooth algebraic curve of genus $g=4$ over $\bbC$ and~$(JC, \Theta_{JC})$ its Jacobian variety. Recall that the theta divisor $\Theta_{JC}$ has precisely two distinct ordinary double points $\pm e$ as singularities.

\begin{lem*} 
There exists a principally polarized abelian scheme $(\calX, \calTheta)$ over a smooth quasi-projective curve $S$,
 \[
 \xymatrix@M=0.4em@C=2.5em@R=2em{
 \calTheta \ar@{^{(}->}[r] \ar@{->>}[dr] & \calX \ar@{->>}[d]^f \\
 & S
 }
 \]
and a point $s\in S(\bbC)$ such that 
\begin{itemize}
 \item the fibre in $s$ is $(\calX_s, \Theta_s) \cong (JC, \Theta_{JC})$, 
 \item the total space $\calX$ and the relative theta divisor $\calTheta$ are nonsingular, 
 \item over $S^* = S \setminus \{s\}$, the structure morphisms $\calX_{S^*} = \calX \times_S {S^*}\rightarrow S^*$ and $\Theta_{\calX^*} = \calTheta \times_S S^* \rightarrow S^*$ are smooth.
\end{itemize}
\end{lem*}

{\em Proof.} Write $JC = \bbC^4 /(\bbZ^4 + \tau_0 \bbZ^4)$ for some point $\tau_0$ in the Siegel upper half plane $\calH_4$, and choose lifts $\pm z_0$ in $\bbC^4$ of $\pm e \in JC$. Let $\pi: \calH_4 \rightarrow \calA_4$ be the analytic quotient map. Since $\pm e$ are ordinary double points of~$\Theta_{JC}$, by the heat equation the gradient of the Riemann theta function~$\theta(\tau, z)$ in $\tau$-direction does not vanish at $(\tau_0, \pm z_0)$.
Since $\calA_4$ is quasi-projective, we can 
use a suitable system of regular parameters of the regular local ring at $\tau_0$ to
construct Zariski-locally a smooth algebraic curve~$S \subset \calA_4$ intersecting $\calJ_4$ transversely in~$s = \pi(\tau_0)$ in a tangent direction along which the gradient of $\tau \mapsto \theta(\tau, \pm z_0)$ does not vanish at~$\tau_0$.  We define $(\calX, \calTheta)$ as the restriction of the universal ppav~$(\bbC^g \times \calH_g)/(\bbZ^{2g} \rtimes \Sp(8, \bbZ)) \rightarrow \calA_4$ to~$S$, shrinking $S$ if required. 
\qed

\subsection{Some notations} \label{subsec:generic_versus_general} In order to transfer the constructions of section~\ref{subsec:theta_square} to the relative situation~\ref{subsec:familyofppavs},
consider the pure Hodge module $\delta = \delta_{\calTheta}[-1] = \bbC_{\calTheta}[3]$.
By Gabber's theorem its relative symmetric resp.~alternating convolution square decomposes as
$S^2(\delta) = \delta_+ \oplus \tau_+$ resp.~$\Lambda^2(\delta) = \one_\calX \oplus \delta_- \oplus \tau_-$
where $\one_\calX$ is a complex supported on the zero section of $\calX$ and concentrated in degree zero, where $\delta_\pm$ are semisimple complexes and
$\tau_{\pm} =  \bigoplus_{\mu  \equiv  \alpha_\pm  (\mathrm{mod}\, 2)}^{} f^*(R^{-3-|\mu|}f_*(\delta_{\calX})) \otimes \delta_\calX[\mu-1]$ as in lemma~\ref{lem:tau}.
For each $t\in S^*(\bbC)$ the restrictions $\tau_\pm |_{\calX_t}$ and $\delta_\pm |_{\calX_t}$ are the complexes that were previously denoted by $\tau_\pm$ resp.~$\delta_\pm$. To prove the main theorem it clearly suffices to show that the new $\delta_\pm \in \Dbc(\calX, \bbC)$ are simple.

\medskip

\label{equivariant}
We consider equivariant perverse sheaves \cite[section III.15]{KW} with respect to the involution $\sigma = -\id_{\calX_s}$. Let $\one_{\sigma \pm}$ be the $\sigma$-equivariant skyscraper sheaf~$\one$ with $\sigma$ acting by $\pm 1$.
For $x\in \calX_s\setminus \{0\}$ let
$\one_{\pm x} = t_x^*(\one) \oplus t_{-x}^*(\one)$
be the simple $\sigma$-equivariant skyscraper sheaf supported in $\{\pm x\}$, with $\sigma$  flipping the two summands.

\subsection{Monodromy filtrations}
\label{subsec:psitheta}
Recall that by construction $\Theta_s= \Theta_{\calX_s}$ is regular except for two distinct ordinary double points~$\pm e$. In particular, then $\delta = \delta_{\calTheta}[-1]$ has unipotent global monodromy by lemma~\ref{Lef}. From section~\ref{MF} we deduce
 that the monodromy filtration diagram of $\Psi_1(\delta)$ is
\[
 \xymatrix@=0.1em{
 &  \; \one_{\pm e}(-2) \;  \ar[dd]^N_\cong \\
   \delta_{\calTheta_s} \; \; \;& \\
 &  \; \one_{\pm e}(-1) \;  \\
 }
\]
since $sp(\delta) = \lambda_{\Theta_s}$ by ex.~\ref{ex:sp} and since the weight filtration of $\lambda_{\Theta_s}$ is defined by the exact sequence 
$0\rightarrow \one_{\pm e} (-1) \rightarrow \lambda_{\Theta_s} \rightarrow \delta_{\Theta_s} \rightarrow 0$
of lemma~\ref{subsec:doublepoints}.

\begin{rem*}
The above implies $\chi(\delta_{\Theta_s}) = \chi(\Psi_1(\delta))-4 = g!-4 = 20$. For any $g\geq 4$ and $r\in \bbN$, similar arguments show that for a ppav $(\calX_s, \Theta_s)$ whose theta divisor $\Theta_s$ has precisely $r$ ordinary double points as singularities, one has $\chi(\delta_{\Theta_s})= g!-2r$. 
\end{rem*}

\label{subsec:psidelta} 

Now consider $\Psi_1(\delta_\pm)$. Clearly $S^2(\one_{\pm e}) = \one_{\pm 2e} \oplus \one_{\sigma+}$ and $\Lambda^2(\one_{\pm e}) = \one_{\sigma-}$, and with notations as in~\cite{WeBN}, theorem 14 on p.~123 in loc.~cit.~and the representation theory of $\Sl(6,\bbC)$ imply that 
\[ S^2(\delta_{\Theta_s}) \;=\; \delta_{3,3} \oplus \delta_{5,1}
 \quad
 \textnormal{and}
 \quad 
 \Lambda^2(\delta_{\Theta_s}) \;=\; \delta_{4,2} \oplus \delta_{6,0}.
\]
In what follows, we write $\delta_\alpha = \pdelt{\alpha} \oplus \cdelt{\alpha}$ where $\cdelt{\alpha}$ is a direct sum of complex shifts of $\delta_{\calX_s}$ and $\pdelt{\alpha}$ has no constituents in $T(X)$. Here $c$ and~$p$ stand for the properties of being {\em constant} resp.~{\em perverse}.

\begin{lem*} 
Up to constituents with vanishing hypercohomology, the complex $\Psi_1(\delta_+)$ has the monodromy filtration diagram
\[ 
 \quad
 \xymatrix@=0.1em{
 & & \\
 & &  (\one_{\pm 2e} \oplus \one_{\sigma+})(-4)   \ar[dd]^N_\cong  \\
 & (\delta_{\Theta_s} * \one_{\pm e})(-2)  \ar[dd]^N_\cong &  \\
  \pdelt{3,3} \oplus \pdelt{5,1} \oplus \one_{\sigma-}(-3) & &  (\one_{\pm 2e} \oplus \one_{\sigma+})(-3)   \ar[dd]^N_\cong  \\
 & (\delta_{\Theta_s} * \one_{\pm e})(-1)  & \\
 & & (\one_{\pm 2e} \oplus \one_{\sigma+}) (-2) \\
 & & 
 }
\]
with graded pieces of weights $4, 5, 6, 7$ and $8$. For $\Psi_1(\delta_-)$ one has the diagram
\[ 
 \xymatrix@=0.1em{
 & & \\
 & &  \one_{\sigma-}(-4) \ar[dd]^N_\cong  \\
 & (\delta_{\Theta_s} * \one_{\pm e})(-2)  \ar[dd]^N_\cong &  \\
  \pdelt{4,2} \oplus (\one_{\pm 2e} \oplus \one_{\sigma+})(-3)  & & \one_{\sigma-}(-3) \ar[dd]^N_\cong  \\
 & (\delta_{\Theta_s} * \one_{\pm e})(-1)  & \\
 & & \one_{\sigma-} (-2)  \\
 & & 
 }
\]
with graded pieces of the same weights as above.
\end{lem*}

{\em Proof.} Use lemma~\ref{lem:GrTensor} and the fact that $\Psi_1(\delta_\pm)$ differs from $\Psi(\delta_\pm)$ at most by constituents with vanishing hypercohomology, taking into account that  $\delta_\pm$ has unipotent global monodromy. To lift the obtained result from $\Pbar(\calX_s)$ to $\Perv(\calX_s)$, note that the complex translates of $\delta_{\calX_0}$ which enter $\Psi_1(S^2(\delta))$ and $\Psi_1(\Lambda^2(\delta))$ can be computed as in lemma~\ref{lem:tau}, so no additional such terms occur in $\Psi_1(\delta_\pm)$.
\qed

\subsection{The degenerate Hodge structure} \label{subsec:Bdegen}
For a suitable choice of $(\calX, \calTheta)$ in lemma~\ref{subsec:familyofppavs} and suitable (fixed) $t\in S^*(\bbC)$ we abbreviate $\Theta_t = \Theta_{\calX_t}$. Then section~\ref{sec:motivic} shows
\[
 \MT(\Theta_t)_{sc} \;=\; \MT(\calX_t) \times \Sp(10, \bbC)
 \quad
 \textnormal{where}
 \quad
 \MT(\calX_t) \;=\; \Sp(8, \bbC).
\]
In section~\ref{app:hypercoh} we compute the natural pure Hodge structure on $\bbH^\bullet(\calX_t, \delta_\pm|_{\calX_t})$ as a representation of this group $\MT(\Theta_t)_{sc}$. But $\bbH^\bullet(\calX_t, \delta_\pm|_{\calX_t})$ can also be equipped with a different (mixed) Hodge structure which is induced from the mixed Hodge modules $\Psi(\delta_\pm)$ via the isomorphism $\bbH^\bullet(\calX_t, \delta_\pm|_{\calX_t})\cong \bbH^\bullet(\calX_s, \Psi(\delta_\pm))$.

\medskip

Let us call this Hodge structure on $\bbH^\bullet(\calX_t, \delta_\pm|_{\calX_t})$ the {\em degenerate} one. Note that the subquotients of the degenerate Hodge structure, such as its invariants under the monodromy operator $N$, are no longer representations of $\MT(\Theta_t)_{sc}$ but only representations of $\MT(\Theta_s)_{sc}=\MT(C)=\Sp(8, \bbC)$.

\medskip

To obtain the degenerate Hodge structure on $\bbH^\bullet(\calX_t, \delta_\pm|_{\calX_t})$, one easily sees that in table~\ref{table:genericHypercoh} of section~\ref{app:hypercoh} the representations of $\MT(\calX_t)=\Sp(8, \bbC)$ are unchanged (including their weights) when viewed as representations of~$\MT(\calX_s)=\Sp(8, \bbC)$. In particular, all of them are invariant under the monodromy operator $N$. However, for the standard representation $B$ of $\Sp(10, \bbC)$ the situation is different as we will see in the lemma below: The degenerate Hodge structure arises from the natural one by pull-back along a monomorphism
\[
 (\id, \varphi): \; \MT(\Theta_s)_{sc} = \Sp(8, \bbC) \hookrightarrow \MT(\Theta_t)_{sc} = \Sp(8, \bbC) \times \Sp(10, \bbC)
\]
for some embedding $\varphi: \Sp(8, \bbC) \hookrightarrow \Sp(10, \bbC)$.

\medskip

Put $\Lambda^i = H^i(\calX_s, \bbQ)$ for $i\in \bbN$. As representations of the group $\Sp(8,\bbC)$ the~$\Lambda^i$ are considered as the exterior powers of the standard representation. Notice that we have $\Lambda^0=(0000), \Lambda^1=(1000)$, $\Lambda^2=(0100)\oplus \Lambda^0$ and $\Lambda^3=(0010)\oplus \Lambda^1$ in the notations of section~\ref{app:hypercoh}.

\begin{lem*} 
On the quotient $B$ of $H^3(\Theta_t, \bbC)\cong  \bbH^0(\calX_t, \delta_{\Theta_t}) \cong \bbH^0(\calX_s, \Psi(\delta))$ with its degenerate limit Hodge structure, the monodromy operator $N$ has the coinvariants 
\[ B/B^N \;=\; \Lambda^0(-2) \] 
which are pure of weight $4$, and the weight filtration of the invariants $B^N$ is given by an exact sequence
\[
 0\longrightarrow \Lambda^0(-1) \longrightarrow B^N \longrightarrow \Lambda^1(-1) \longrightarrow 0
\]
with $\Lambda^1(-1)$ pure of weight $3$ and with $\Lambda^0(-1)$ pure of weight $2$. 
\end{lem*}

{\em Proof.}  Since $\Phi(\delta)=\one_{\pm e}(-2)$ and $sp(\delta) = \lambda_{\Theta_s}$ by section~\ref{subsec:psitheta}, we get from proposition~\ref{cor:sp}(a) an exact sequence
\[
 0 \longrightarrow \bbH^0(\calX_s, \lambda_{\Theta_s})
 \stackrel{\alpha}{\longrightarrow}  \underbrace{\bbH^0(\calX_t, \delta_{\Theta_t})}_{= \Lambda^3 \oplus B}
 \stackrel{\beta}{\longrightarrow} \underbrace{\bbH^0(\calX_0, \one_{\pm e}(-2))}_{=(0000)^{\oplus 2}}
\]
where $\alpha$ and $\beta$ are morphisms of Hodge structures if the middle term is equipped with the degenerate limit Hodge structure.
By the local invariant cycle theorem the image of $\alpha$ is the subspace of $N$-invariant elements. Since $\dim_\bbC (B)=10$ and $\dim_\bbC(\Lambda^0) + \dim_\bbC (\Lambda^1) = 9$, by a dimension count the claim follows once we can exhibit an exact sequence
\begin{equation} \label{eq:degenB}
 0 \longrightarrow \Lambda^0(-1) \longrightarrow \bbH^0(\calX_s, \lambda_{\Theta_s})
 \longrightarrow \Lambda^3 \oplus \Lambda^1(-1) \longrightarrow 0
\end{equation}
The sequence~(\ref{eq:degenB}) is obtained as the short exact sequence
\begin{equation} \label{eq:ses2}
0 \longrightarrow \ker(\bbH^0(\nu)) \longrightarrow \bbH^0(\calX_s, \lambda_{\Theta_s}) \longrightarrow \bbH^0(\calX_s, \delta_{\Theta_s}) \longrightarrow 0
\end{equation}
obtained by splicing up the long exact sequence attached to the short exact sequence $0 \to \one_{\pm e}(-1) \to \lambda_{\Theta_s} \to \delta_{\Theta_s} \to 0$ of lemma~\ref{subsec:doublepoints}, i.e.~with $ \ker(\bbH^0(\nu)) $ isomorphic to the quotient
\begin{equation} \label{eq:ses3}
 0 \longrightarrow \bbH^{-1}(\calX_s, \lambda_{\Theta_s}) \longrightarrow \bbH^{-1}(\calX_s, \delta_{\Theta_s}) \longrightarrow  \ker(\bbH^0(\nu)) \to 0.
\end{equation}
The right hand side of (\ref{eq:ses2}) and the middle term of~(\ref{eq:ses3}) are given by the cohomology of the smooth curve $C$ via $\bbH^\bullet(\delta_{\Theta_s}) = \Lambda^{g-1}(\bbH^\bullet(\delta_C))$, see \cite{WeBN}. Thus
\[
 \bbH^0(\calX_s, \delta_{\Theta_s}) \;=\; \Lambda^3 \oplus \Lambda^1(-1)
 \quad
 \textnormal{and}
 \quad
 \bbH^{-1}(\calX_s, \delta_{\Theta_s}) \;=\; \Lambda^2 \oplus \Lambda^0(-1).
\]
To finish the proof it remains to compute $\bbH^{-1}(\calX_s, \lambda_{\Theta_s})$ in~(\ref{eq:ses3}). For this use the long exact cohomology sequence associated with the triple $(\calX_s \setminus \Theta_s, \calX_s, \Theta_s)$. Since we have $H_c^i(\calX_s\setminus \Theta_s, \bbC)=0$ for $i\leq 3$ by the ampleness of the theta divisor, we get
\[
 \bbH^{-2}(\calX_s, \lambda_{\calX_s}) \; = \; H^2(\calX_s, \bbC) \; \stackrel{\sim}{\longrightarrow} \;
 H^2(\Theta_s, \bbC) \;=\; \bbH^{-1}(\calX_s, \lambda_{\Theta_s})
\] 
which is isomorphic to $\Lambda^2$. From this one easily concludes the proof.
\qed

\subsection{Two big constituents} \label{subsec:bookkeeping1} 
We now exhibit two irreducible constituents~$\gamma_\pm \hookrightarrow \delta_\pm$ such that~$\spbar(\gamma_\pm)$ differs from $\spbar(\delta_\pm)$ at most by skyscraper sheaves. For $w_0\in \bbZ$ and a mixed perverse sheaf $\pi$ on $\calX_s$, we will denote by $\pi_{w<w_0}$ the maximal perverse subsheaf of $\pi$ of weights $< w_0$.

\begin{keylem*} \label{prop:big}
There exist unique irreducible constituents $\gamma_\pm \hookrightarrow \delta_\pm$ in $\Perv(\calX_t)$ with the property
$\pdelt{5,1} \oplus \pdelt{3,3} \hookrightarrow \spbar(\gamma_+)$
and
$\pdelt{4,2} \hookrightarrow \spbar(\gamma_-)$ in $\Perv(\calX_s)$.
The perverse sheaves $sp(\gamma_+)_{w<6}$ and $sp(\gamma_-)_{w<6}$ both admit $\delta_{\Theta_s}*\one_{\pm e}(-1)$ as a quotient.
\end{keylem*}

{\em Proof.}
Recall from lemma~\ref{subsec:psitheta} that $K=\delta_{\Theta_s} * \one_{\pm e}(-1)$ is a $\sigma$-equivariant simple constituent of $sp(\delta_\pm)_{w<6}$. Suppose $\epsilon_\pm$ is a constituent of $\delta_\pm$ for which $sp(\epsilon_\pm)$ does not contain $K$. Then $sp(\epsilon_\pm)_{w<6}$ and $\Phi(\epsilon_\pm)$ are skyscraper sheaves, hence
\[
 \bbH^{-3}(\calX_s, \Phi_1(\epsilon_\pm)) \; = \; 0
 \quad \textnormal{and} \quad
 \bbH^{-2}(\calX_s, sp(\epsilon_\pm)) \; \stackrel{\sim}{\longrightarrow} \;
 \bbH^{-2}(\calX_s, \Psi_1(\epsilon_\pm)).
\]
By the local invariant cycle theorem we then have an isomorphism
\[
 \bbH^{-2}(\calX_s, \Psi_1(\epsilon_\pm)) \; \stackrel{\sim}{\longrightarrow} \;
 \bbH^{-2}(\calX_t, \epsilon_\pm |_{\calX_t})^N.
\]
To compute $\bbH^{-2}(\calX_t, \epsilon_\pm |_{\calX_t})^N$ we use the second line of table~\ref{table:genericHypercoh} in section~\ref{app:hypercoh} where $\bbH^{-2}(\calX_t, \delta_\pm |_{\calX_t})$ is listed. The monodromy operator $N$ acts non-trivially only on $B$, and $B/B^N$ and $B^N$ were computed in lemma~\ref{subsec:Bdegen}. 

\medskip

We claim that the summand $(1000)\otimes B^N$ of $\bbH^{-2}(\calX_t, \delta_\pm |_{\calX_t})^N$ is linearly disjoint from $\bbH^{-2}(\calX_t, \epsilon_\pm |_{\calX_t})^N$. Otherwise the summand $(1000)\otimes B^N$ would be contained in $\bbH^{-2}(\calX_s, \spbar(\epsilon_\pm))$ by global monodromy reasons since $\epsilon_\pm$ are complexes defined globally on $\calX$. Namely, $(1000)\otimes B$ would occur in $\bbH^\bullet(\calX_s, \Psi(\epsilon_\pm)) = \bbH^\bullet(\calX_t, \epsilon_\pm)$ as a representation of the Mumford-Tate group, so by the local invariant cycle theorem the representation $(1000)\otimes B^N$ would occur in $\bbH^{-2}(\calX_s, \spbar(\epsilon_\pm))$. This is impossible, since $B^N$ is not pure.

\medskip

By section~\ref{app:hypercoh}, our claim shows
$\bbH^{-2}(\calX_s, \spbar(\epsilon_\pm)) \subseteq (1010) \oplus (0100)$. 
Hence again by section~\ref{app:hypercoh}, the simple perverse sheaves
\[
 L \;=\; \pdelt{5,1}, \; \pdelt{4,2} \quad \textnormal{or} \quad \pdelt{3,3}
\]
cannot appear as constituents of $\spbar(\epsilon_\pm)$ because for these perverse sheaves $L$ the representation of the group $\MT(\calX_s)=\Sp(8, \bbC)$ on $\bbH^{-2}(\calX_s, L)$ is not contained in $(1010) \oplus (0100)$ (see line $2$ of table~\ref{table:specialHypercoh}).

\medskip

Hence for any $\sigma$-equivariant constituent $\gamma_\pm$ of $\delta_\pm$ such that $\spbar(\gamma_\pm)$ contains one of the constituents $L$ above, $sp(\gamma_\pm)$ has $K$ as a constituent. Since from lemma~\ref{subsec:psitheta} we know $K$ enters with multiplicity one in $sp(\delta_\pm)$, it follows that there are unique $\sigma$-equivariant simple constituents $\gamma_\pm$ of $\delta_\pm$ containing one of the simple perverse sheaves $L$ above. These satisfy
\[
 \pdelt{5,1}\oplus \pdelt{3,3} \; \hookrightarrow \; \spbar(\gamma_+)
 \quad \textnormal{and} \quad
 \pdelt{4,2} \; \hookrightarrow \; \spbar(\gamma_-),
\]
which easily concludes the proof.
\qed

\subsection{How to exclude skyscraper sheaves}
\label{subsec:no-skyscraper}
Now define $\varepsilon_\pm \in \Dbc(\calX, \bbC)$ by the decomposition $\delta_\pm = \gamma_\pm \oplus \varepsilon_\pm$ of semisimple perverse sheaves.
To prove the main theorem we must show $\varepsilon_\pm =0$.
From lemmas~\ref{prop:big} and~\ref{subsec:psidelta} we know 
\[
\begin{array}{rclrcl}
 sp_{w<6}(\varepsilon_+) &\hookrightarrow& (\one_{\pm 2e}\oplus \one_{\sigma+})(-2), 
 & \spbar(\varepsilon_+) &\hookrightarrow& \one_{\sigma-}(-3), \\ 
 sp_{w<6}(\varepsilon_-) &\hookrightarrow& \one_{\sigma-}(-2),  
 &\spbar(\varepsilon_-) &\hookrightarrow& (\one_{\pm 2e}\oplus \one_{\sigma+})(-3).
\end{array}
\]
In particular, $sp(\varepsilon_\pm)$ and hence $\Psi(\varepsilon_\pm)$ are skyscraper sheaves, so $\bbH^0(\calX_s, \varepsilon_\pm)$ is a direct sum of trivial representations $(0000)$. A look at table~\ref{table:genericHypercoh} shows that these can only arise from lines $2$ and $4$ of this table; possible candidates arising from $(1000)\otimes B$ in line $2$ are ruled out by global monodromy reasons since $B$ is an irreducible representation of $\MT(\Theta_t)$. The irreducibility of $\one_{\pm 2e}$ as a $\sigma$-equivariant perverse sheaf hence implies
\[
 sp(\varepsilon_+) \;=\; sp_{w<6}(\varepsilon_+) \;\hookrightarrow\; \one_{\sigma+}(-2)
 \quad
 \textnormal{and}
 \quad
 sp(\varepsilon_-) \;=\; \spbar(\varepsilon_-) \; \hookrightarrow\; \one_{\sigma+}(-3).
\]
We now claim $sp(\varepsilon_+) = 0$. Indeed, otherwise $\bbH^0(\calX_s, sp(\varepsilon))$ would be the trivial representation $(0000)$ of weight $4$. But the trivial representation could only arise from one of the summands
\[
 (S^2(B))^N, \quad (0010)\otimes B^N \quad \textnormal{or} \quad (0000)
\]
in line $4$ of the first column of table~\ref{table:genericHypercoh}. The first two of these summands cannot contribute because of global monodromy reasons, and the last one has the wrong weight $6$. This proves our claim that $sp(\varepsilon_+)=0$. But then also $\Psi(\epsilon_+)=0$ and hence
\[
 \bbH^\bullet(\calX_t, \varepsilon_+|_{\calX_t}) \;=\; \bbH^\bullet(\calX_s, \Psi(\varepsilon_+)) \;=\; 0
\]
which easily implies $\varepsilon_+=0$ since $\varepsilon_+|_{\calX_t} \notin T(\calX_t)$.

\medskip

By the same argument, to show $\varepsilon_-=0$ it suffices to see that $sp(\varepsilon_-)=0$. If $sp(\varepsilon_-)\neq 0$, then by what we have seen above
\[
 sp(\varepsilon_-) \;=\; \spbar(\varepsilon_-) \;=\; \one_{\sigma+}(-3),
\]
hence also $\Psi(\varepsilon_-)=\one_{\sigma+}(-3)$ and therefore $\bbH^\nu(\calX_t, \varepsilon_-|_{\calX_t})$ is $\bbC$ for $\nu = 0$ and zero otherwise. Therefore $\chi(\varepsilon_-|_{\calX_t})=1$, and by~\cite{KrW2} then $\varepsilon_-|_{\calX_t} = \one$. But this is impossible since $\delta_{\Theta_t}*\delta_{\Theta_t}$ contains the unit object $\one$ only with multiplicity one.
\qed

\section{Equivalence of conjectures~\ref{conj:thetasquare} and \ref{conj:thetagroup}} \label{sec:super-mackey}

In this section we prove our earlier statement that the conjectures~\ref{conj:thetasquare} and~\ref{conj:thetagroup} are equivalent for any~$g$. Again let $k=\Qbar_l$ or $k=\bbC$.\label{subsec:mackey} 
Let $G$ be a reductive group over $k$ and let $H\hookrightarrow G$ be a closed subgroup of finite index. For a representation $U$ of $G$ we will denote by $R^G_H(U)$ the restriction of $U$ to $H$. Similarly, for a representation~$V$ of~$H$, let $I_H^G$ denote the induced representation of $G$. One then easily proves the following version of Mackey's lemma.

\begin{lem*}
For every irreducible representation $U$ of $G$ there is a subgroup $H'\subseteq G$ containing $H$ and an irreducible representation $V'$ of $H'$ with the property that the restriction~$R_H^{H'}(V')$ is isotypic and $U \cong  I_{H'}^G(V')$.
\end{lem*}

\begin{cor*}
If $U$ is an irreducible representation of $G$ and $R^G_{G^0}(U)$ contains the trivial representation, then all constituents of $R^G_{G^0}(U)$ are trivial.
\end{cor*}

Using this we now prove the equivalence of conjectures~\ref{conj:thetasquare} and~\ref{conj:thetagroup} for all~$g$. \label{equ}
We must see that~\ref{conj:thetasquare} implies~\ref{conj:thetagroup}. If~\ref{conj:thetasquare} holds, the irreducible representations $U=\omega(\delta_\Theta)$ and $W_\pm=\omega(\delta_\pm)$ of $G=G(X, \Theta)$ satisfy 
\[
 S^2(U) \;=\; \begin{cases}
               W_+ \\
		W_+ \oplus k
              \end{cases}
\textnormal{and}
\;\;
\Lambda^2(U) \;=\; \begin{cases}
                    W_- \oplus k \\
			W_-
		\end{cases}
\textnormal{for}
\;\;
\begin{cases}
 \textnormal{$g$ even}, \\
 \textnormal{$g$ odd}.
\end{cases}
\]
Then $U$ cannot be a representation induced from a proper subgroup of $G$ because otherwise $S^2(U)$ would contain at least two non-trivial irreducible constituents. Hence by the isotypic case of the lemma,
$R^G_{G^0}(U) = U_0^{\oplus n}$
for some irreducible super representation $U_0$ of $G^0$ and some $n\in \bbN$. 

\medskip

If $n>1$, then for $g$ even resp.~odd, $R^G_{G^0}(W_+)$ resp.~$R^G_{G^0}(W_-)$ contains $U_0\otimes U_0$. Since by self-duality of $U_0$ the trivial representation enters $U_0\otimes U_0$ and since $W_+$ and $W_-$ are irreducible, the corollary would then imply that $U_0$ were trivial.
By faithfulness then $G^0$ would be trivial, i.e.~$G$ were a finite group. Then we could find $K\in \Perv(X)$ with $\omega(K)$ being the regular representation $R=k[G]$. Since $R\otimes R = R^m$ for $m= |G|$, then $K*K \equiv K^{\oplus m}$ modulo $T(X)$, so $\bbH^i(X, K)=0$ for all~$i\neq 0$ by~\cite[lemma 5, p.~17]{WeBN}, a contradiction since $\delta_\Theta$ is a direct summand of~$K$ (every representation of $G$ enters the regular representation $R$).

\medskip

Hence $n=1$ and $R^G_{G^0}(U)=U_0$ is irreducible. Then $R^G_{G^0}(U)$ contains a unique highest weight vector up to scalars, so the same holds for $R^G_{G^0}(S^2(U))$. Thus $R^G_{G^0}(W_+)$ is neither induced nor isotypic in a non-trivial way, hence by the lemma it must be irreducible.
It follows that $W_+$ is irreducible as a super representation of the super Lie algebra $\frakg$ of $G$. The classification in~\cite{KrW} and the fact that $\dim(V)=g!$ then imply that
\begin{enumerate}
 \item if $g$ is even, $G$ is of type $A_{g!-1}$ or of type $C_{g!/2}$,
 \item if $g$ is odd, $G$ is of type $D_{g!/2}$,
\end{enumerate}
and that in all these cases $V$ is the standard representation.
On the other hand, we have already observed in section \ref{subsec:thetagroupmotivation} that the faithful action of $G$ on $V$ preserves a nondegenerate alternating resp.~symmetric bilinear form $\beta: V\times V \rightarrow k$, defining an embedding of $G$ into $\Sp(V, \beta)$ resp.~$\mathrm{O}(V, \beta)$, for $g$ even resp.~odd. So $G=\Sp(V, \beta)$ in case (a). In case (b) either $G=\mathrm{O}(V, \beta)$ or $G=\SO(V, \beta)$; but the arguments of~\cite{KrW2} show that $G$ does not admit any non-trivial character, so $G=\SO(V, \beta)$ and we are done.
\qed

\section{Appendix: Two lemmas on perverse sheaves} \label{app:perv}

Let $k=\Qbar_l$ or $k=\bbC$, and for a variety $Y$ over an algebraically closed field of characteristic zero, let $k_Y\in \Dbc(Y, k)$ denote the constant sheaf on $Y$. 

\subsection{Weight filtrations} For the computation of nearby cycles the following finer version of \cite[lemma~2]{WeBN} has been used.

\begin{lem*} \label{subsec:doublepoints}
(i) If $Y$ is an irreducible normal local complete intersection of dimension~$d$ with singular locus $\Sigma \subset Y$, then on $Y$ we have an exact sequence of perverse sheaves
\[
 0 \longrightarrow \psi_Y \longrightarrow k_Y[d] \longrightarrow \delta_Y \longrightarrow 0
\]
where $\psi_Y$ is a mixed perverse sheaf of weights $<d$ whose support is contained in~$\Sigma$ and $\delta_Y=\IC_Y[d]$ is a pure perverse sheaf of weight $d$.

\medskip

(ii) If furthermore the only singularities of $Y$ are finitely many ordinary double points $y_1, \dots, y_n$, then 
$
 \psi_Y = \oplus_{i=1}^n \delta_{y_i}(-\frac{d-1}{2})$
is a direct sum of skyscraper sheaves for $d$ odd and $\psi_Y = 0$ otherwise.
\end{lem*}

{\em Proof.} {\em (i)} This follows from lemma 1 and 2 in~\cite{WeBN}.  For the convenience of the reader we sketch the proof.  Suppose $Y$ is irreducible and normal. Then as in loc. cit. there exists a natural morphism 
$\nu: \lambda_Y \to \delta_Y$ of sheaf complexes, such that $\calH^{-d}(\nu)$ is an isomorphism.
For $Y$ normal $\lambda_Y = k_Y[d]$, hence $\lambda_X$ is a complex of weights $\leq d$ whereas $\delta_Y=\IC_Y[d]$ is pure of weight $d$. Since $\lambda_Y \in {}^p D^{\leq 0}(Y)$ by definition and $\delta_Y$ is a perverse sheaf, $\nu$ factorizes over the truncation morphism 
$\lambda_Y \to {}^p H^0(\lambda_Y)$. 
Since $\delta_Y$ is an irreducible perverse sheaf and $\nu$ is nontrivial, it is easy to see that the induced  morphism $\mu: {}^p H^0(\lambda_Y) \to \delta_Y$ is nontrivial.
Hence $\mu$ defines an epimorphism in the category of perverse sheaves. So we obtain a distinguished triangle
\[
0 \longrightarrow \psi_Y \longrightarrow k_Y[d] \longrightarrow \delta_Y \longrightarrow 0
\ . \]
with $\psi_Y\in {}^p D^{\leq 0}(Y)$.  The long exact sequence of cohomology sheaves
for this distinguished triangle implies $\calH^{-\nu}(\psi_Y) \cong \calH^{-\nu -1}(\delta_Y)$
for $\nu \geq -d+2$ and $\calH^{-\nu}(\psi_Y) =0$ otherwise, and 
hence $\psi_Y$ is
of weights $<d$. Finally,
if $Y$ is also a local complete intersection, it follows from~\cite[III.6.5]{KW} that $\lambda_Y$
itself is a perverse sheaf.

\medskip

{\em (ii)} Let $\pi: \Ytilde \rightarrow Y$ be the blow-up of $Y$ in $\Sigma = \{y_1, \dots, y_n\}$. Notice that $\Ytilde$ is smooth. Since $\pi$ restricts to an isomorphism over $U=Y\setminus \Sigma$, it follows from purity that
$R\pi_*(\lambda_\Ytilde) = \delta_Y \oplus \gamma$
for some $\gamma \in \Dbc(Y, k)$
with $\mathit{Supp}(\gamma)\subseteq \Sigma$. Let us write
$
 \gamma = \bigoplus_{i=1}^n \bigoplus_{j\in \bbZ}  (\delta_{y_i}((j-d)/2)[j] )^{\oplus m_{ij}}$
with $m_{ij}\in \bbN_0$.
Since the $y_i$ are ordinary double points, the $Q_i = \pi^{-1}(y_i)$ are smooth quadrics of dimension $d-1$, so by~\cite[exp.~XII, th.~3.3]{SGA7}
\[ H^j(Q_i, k) \;=\; 
 \begin{cases}
  k (-\frac{j}{2}) & \textnormal{for $j \in \{ 0,2,\dots, 2(d-1)\} \setminus \{ d-1\}$}, \\
  \bigl(k (-\frac{j}{2})\bigr)^{2\delta} & \textnormal{for $j=d-1$}, \\
  0 & \textnormal{otherwise}, 
 \end{cases}
\]
where $\delta = 1$ for $d$ odd and $\delta = 0$ for $d$ even. 
Now 
$\calH^\bullet(R\pi_*(\lambda_\Ytilde))_{y_i}  =  H^\bullet(Q_i, k)[d]$.
In particular, $\calH^d(\gamma)_{y_i}=0$ and $\calH^{d-2}(\gamma)_{y_i} = k$, so~$m_{i,d}=0$ and~$m_{i,d-2}=1$.
Then $m_{i,\pm d}=0$ and $m_{ij}\geq 1$ for $j=d-2, d-4, \dots, 4-d, 2-d$ by the hard Lefschetz theorem.
Another look at stalk cohomology then shows that $m_{ij}=1$ for $j=d-2, d-4, \dots, 4-d, 2-d$ and $m_{ij}=0$ otherwise. Thus
\[
 \calH^r(\delta_Y)_y \;=\; 
 \begin{cases}
  k & \textnormal{for $r=-d$ and all $y$},\\
 \bigl( k (-\frac{d-1}{2}) \bigr)^{\delta}& \textnormal{for $r=-1$ and $y\in \Sigma$},\\
  0 & \textnormal{otherwise},
 \end{cases}
\]
with $\delta$ as above.
On the other hand, $k_Y[d] = \lambda_Y = \calH^{-d}(\delta_Y) = \sttau_{\leq -d}(\delta_Y)$
by normality of $Y$ \cite[lemma 1]{WeBN} and \cite[III.5.14]{KW}. 
\qed

\subsection{Counting IC-constituents}\label{subsec:lowest-hypercohomology}
Let $Y$ be a variety of dimension $g$ over an algebraically closed field of characteristic zero, and let $\delta \in \Perv(Y)$.

\begin{lem*}
The perverse sheaf $\delta$ admits~$\delta_Y$ as a constituent iff $\bbH^{-g}(Y, \delta)\neq 0$.
\end{lem*}

{\em Proof.} We know that $\delta \in \stD^{\geq -g}(Y, \bbC)$, so $E_2^{p,q}\!=\! H^p(Y, \calH^q(\delta))\! \Rightarrow \!\bbH^{p+q}(Y, \delta)$ implies that $\bbH^{-g}(Y, \delta)=H^0(Y, \calH^{-g}(\delta))$. Let
$j: U\hookrightarrow Y$ be an open dense subset such that~$\calG = \delta[-g]|_U$ is locally constant, possibly zero. Then~\cite[III.5.14]{KW} shows that $\calH^{-g}(\delta)=j_*( \calG)$, hence $H^0(Y, \calH^{-g}(\delta))=H^0(U, \calG)$, and this group is zero if and only if $\calG$ has no constant subsheaf. 
\hfill \qed \\

\section{Appendix: Hypercohomology computations} \label{app:hypercoh}

In this appendix we determine the hypercohomology of some perverse sheaves required in section~\ref{sec:decomposition}. Let $(X, \Theta)$ be a general ppav in~$\calA_4$ as in section~\ref{sec:motivic}, and define $\delta_\pm \in \Perv(X)$ as in~\ref{subsec:theta_square}. Let $JC$ be the Jacobian of a general curve $C$ of genus $g=4$, and consider the associated perverse sheaves~$\pdelt{\alpha}$ on $JC$ as in section \ref{subsec:psidelta}.
Let us group the hypercohomology of perverse sheaves into packages $[n]_t$ which are stable under the Lefschetz operator and occur precisely in degrees $n, n-2, \dots, 2-n, -n$ for some $n\in \bbN_0$.  Denote by $(a_1,a_2,a_3,a_4)$ the irreducible representation of $\Sp(8,\bbC)$ with highest weight~$a_1\omega_1 + \cdots + a_4\omega_4$ for the fundamental dominant weights $\omega_1, \dots, \omega_4$, and let $B$ be as in section~\ref{sec:motivic} for $(X, \Theta)$. To indicate that a representation enters with multiplicity $m > 1$ we use a superscript ${\oplus m}$, and we specify the action of $\sigma=-\id_{X}$ with a a subscript $\sigma \pm$. 

\begin{table}
\[ \footnotesize
\begin{array}{|c||l|l|} \hline \medrowheight
 \bullet & \bbH^\bullet(X, \delta_+) & \bbH^\bullet(X, \delta_-) 
\\ \hline 
 [3]_t & (0000)_{\sigma +}\otimes B & (0000)_{\sigma +}\otimes B  \\ 
 &  (0010)_{\sigma -} &  \\ \hline 
 [2]_t &  (1000)_{\sigma -}\otimes B & (1000)_{\sigma -}\otimes B  \\
 &  (1010)_{\sigma +} & (0100)_{\sigma +}  \\ 
 &  (0100)_{\sigma +} & \\ \hline
 [1]_t & (0100)_{\sigma +}\otimes B & (0100)_{\sigma +}\otimes B \\
 & (1100)_{\sigma -} & (1100)_{\sigma -} \\ 
 & (1000)_{\sigma -}^{\oplus 2} & (1000)_{\sigma -} \\
 & (0110)_{\sigma -} &  \\ \hline 
 [0]_t & S^2(B) & \Lambda^2(B)  \\
 &  (0010)_{\sigma -} \otimes B &  (0010)_{\sigma -}\otimes B  \\
 & (2000)_{\sigma +}^{\oplus 2} & (2000)_{\sigma +} \\
 &  (0200)_{\sigma +} & (0200)_{\sigma +}  \\
 &  (0020)_{\sigma +} & (0000)_{\sigma +}  \\ 
 &  (0000)_{\sigma +} &  \\ \hline
\end{array}
\]
\caption{Decomposition of the Hodge structures $\bbH^\bullet(X, \delta_\pm)$ as representations of $\MT(\Theta)_{sc} = \Sp(8, \bbC) \times \Sp(10, \bbC)$} \label{table:genericHypercoh}
\end{table}

\begin{table}
\[ \footnotesize
\begin{array}{|c||l|l|l|} \hline \medrowheight
 \bullet & \bbH^\bullet(JC, \pdelt{5,1}) & \bbH^\bullet(JC, \pdelt{4,2}) & \bbH^\bullet(JC, \pdelt{3,3}) 
\\ \hline 
 [3]_t &  & (1000) & (1000) \\ 
 &   & & (0010) \\ \hline 
 [2]_t &   (0000) & (2000) & (2000) \\
 &   & (0100)^{\oplus 2} & (1010) \\ 
 &   & (0000)^{\oplus 3} & (0100)^{\oplus 2} \\
 & & & (0000) \\ \hline 
 [1]_t &  (1000)^{\oplus 2} & (1100)^{\oplus 2} & (1100)^{\oplus 2}\\
 &  & (1000)^{\oplus 4} & (1000)^{\oplus 3} \\ 
 &  & (0010) & (0110) \\
 &  & & (0010) \\ \hline 
 [0]_t &  (2000) & (2000) & (2000)^{\oplus 2} \\
 &   (0100) & (1010) & (1010) \\
 & (0000) & (0200) & (0200) \\
 &  & (0100)^{\oplus 3} & (0100) \\
 &   & (0000)^{\oplus 2} & (0020) \\ 
 &  & & (0000)^{\oplus 2} \\ \hline
\end{array}
\]
\caption{Decomposition of the Hodge structure $\bbH^\bullet(JC, \pdelt{\alpha})$ as representations of $\MT (C) = \Sp (8, \bbC)$}
\label{table:specialHypercoh}
\end{table}

\begin{lem*} \label{lem:motivic2}
The representation of $\MT(\Theta)_{sc} = \Sp(8, \bbC)\times \Sp(10, \bbC)$ on the Hodge structure $\bbH^\bullet(X, \delta_\pm)$ is given in table~\ref{table:genericHypercoh}, with $\sigma = -\id_X$ acting trivially on $B$.
The representation of $\MT(C)=\Sp(8, \bbC)$ on $\bbH^\bullet(JC, \pdelt{\alpha})$ for $\alpha \in \{ (5,1), (4,2), (3,3)\}$ is given in table~\ref{table:specialHypercoh}.
\end{lem*}

{\em Proof.}
This has been worked out using the computer algebra systems {\sffamily MAGMA} and {\sffamily SAGE}.
For the case of $\bbH^\bullet(X, \delta_\pm )$ take the symmetric resp.~alternating square of~$
 \bbH^\bullet(X, \delta_{\Theta}) =
 (0000)[3]_t \oplus (1000)[2]_t \oplus (0100)[1]_t \oplus ((0010)\oplus B))[0]_t
$
in the super sense 
and then subtract $\bbH^\bullet(X, \tau_\pm)$. Here  $\sigma$ acts by $-1$ on~$(1000)$ and on~$(0010)$ but trivially on~$(0000)$, $(0100)$, $(0001)$. To check that $\sigma$ acts trivially on $B$, note that by \cite[ex.~4.12(14)]{BL} the number of $2$-torsion points on $\Theta$ is $2^{g-1}(2^g - 1) = 120$; the Lefschetz fixed point formula for~$\sigma$ then implies that $\sigma^*|_B = \id_B$.
For the last three columns, for $a\geq b > 0$ the Littlewood-Richardson rule in~\cite{WeBN} says that~$\pdelt{a,b}$ is the difference of~$\pdelt{a}*\pdelt{b}$ and $\pdelt{a+1}*\pdelt{b-1}$ up to complex shifts of~$\delta_{JC}$. These complex shifts of $\delta_{JC}$ can be computed as in the proof of lemma~\ref{lem:tau}.
\qed

\end{document}